\documentclass[11pt]{article}
\usepackage{amsmath}
\usepackage{amsthm}
\usepackage{amssymb}
\usepackage{graphicx}
\usepackage{hyperref}
\usepackage{amsfonts}
\usepackage{cite}
\usepackage{color}

\numberwithin{equation}{section}

\newtheorem{theorem}{Theorem}[section]
\newtheorem{lemma}{Lemma}[section]

\newtheorem{corollary}{Corollary}[section]

\newtheorem{remark}{Remark}

\newcommand{\fd}[1]{\mathcal{D}^{#1}_t}
\newcommand{\dfd}[1]{\mathcal{D}^{#1}_\tau}

\newcommand{\diff}{\triangledown_\tau}
\newcommand{\taumax}{\tau}

\newcommand{\opL}{\mathcal{L}}

\newcommand{\Gloc}{G_{\mathrm{loc}}}
\newcommand{\Ghis}{G_{\mathrm{his}}}

\newcommand{\Ass}[1]{\textbf{\upshape A#1}}
\newcommand{\Mss}[1]{\textbf{\upshape M#1}}
\newcommand{\defeq}{:=}
\newcommand{\zd}{\,\mathrm{d}}

\newcommand{\abs}[1]{\left|#1\right|}
\newcommand{\absb}[1]{\big|#1\big|}

\newcommand{\abst}[1]{|#1|}
\newcommand{\bra}[1]{\left(#1\right)}
\newcommand{\brab}[1]{\big(#1\big)}
\newcommand{\braB}[1]{\Big(#1\Big)}

\newcommand{\kbrab}[1]{\big[#1\big]}
\newcommand{\kbraB}[1]{\Big[#1\Big]}

\newcommand{\iprod}[1]{\left\langle#1\right\rangle}

\newcommand{\myinnerb}[1]{\big\langle#1\big\rangle}

\newcommand{\mynorm}[1]{\left\|#1\right\|}
\newcommand{\mynormb}[1]{\big\|#1\big\|}

\title{Sharp $H^1$-norm error estimates of two time-stepping \\schemes
for reaction-subdiffusion problems}
\author{Jincheng Ren\thanks{College of Mathematics and Information Science, Henan University
of Economics and Law, Zhengzhou 450046, P. R. China. E-mail:
renjincheng2001@126.com.  This author is supported in part by NSFC
grant 11601119, the program No.18HASTIT027 for HASTIT, and Young
talents Fund of HUEL.}
\quad Hong-lin Liao\thanks{Corresponding
author. Department of Mathematics, Nanjing University of Aeronautics
and Astronautics,  Nanjing 211106, P. R. China. E-mail:
liaohl@csrc.ac.cn.  This
 author is supported in part by the grant
1008-56SYAH18037 from NUAA Scientific Research Starting Fund of
Introduced Talent, and a grant DRA2015518 from 333 High-level
Personal Training Project of Jiangsu Province.}
\quad Jiwei
Zhang\thanks{Beijing Computational Science Research Center, Beijing
100193, P. R. China. E-mail: jwzhang@csrc.ac.cn.  This author is
supported in part by NSFC grants 11771035, 91430216 and NSAF U1530401. }
\quad Zhimin Zhang\thanks{Beijing Computational Science Research Center, Beijing 100193, P. R. China;
and Department of Mathematics, Wayne State University, Detroit, MI 48202, USA.
E-mails: zmzhang@csrc.ac.cn  and ag7761@wayne.edu.  This
author is supported in part by NSFC grants 11471031, 91430216, U1530401,
and by NSF under grant DMS-1419040.}}
\date{\today}

\begin{document}

\maketitle

\begin{abstract}
Due to the intrinsically initial singularity of solution and the discrete convolution form in numerical Caputo derivatives,
the traditional $H^1$-norm analysis (corresponding to the case for a classical diffusion equation) to the time approximations
of a fractional subdiffusion problem always leads to suboptimal error estimates (a loss of time accuracy).
To recover the theoretical accuracy in time, we propose an improved discrete Gr\"{o}nwall inequality and apply it
to the well-known L1 formula and a fractional Crank-Nicolson scheme. With the help of
a time-space
error-splitting technique and the global consistency analysis,
 sharp $H^1$-norm error estimates of the two nonuniform approaches are established for a reaction-subdiffusion problems.
Numerical experiments are included to confirm the sharpness of our analysis.
\vspace{1em}\vspace{1em}
\end{abstract}

\section{Introduction}
\setcounter{equation}{0} Sharp $H^1$-norm error estimates are
established for two nonuniform time approximations to a linear
reaction-subdiffusion problems \cite{Hilfer:2000,Podlubny:1999} in a spatial domain~$\Omega$
\begin{equation}\label{eq: IBVP}
\begin{aligned}
\fd{\alpha} u+\opL u&=c(x) u +f(x, t)&&\text{for $x\in\Omega$ and $0<t\le T$,}\\
u&=u_b(t,x)&&\text{for $x\in\partial\Omega$ and $0<t<T$,}\\
u&=u_0(x)&&\text{for $x\in\Omega$ when $t=0$,}
\end{aligned}
\end{equation}
where $\opL$ is a linear, second-order, strongly-elliptic partial
differential operator in the spatial variable~$x$, and $c(x)$ is a
reaction coefficient satisfying $\abs{c(x)}\le \kappa$ for a
positive constant $\kappa$. Here,
$\fd{\alpha}={}^{\text{C}}_0\fd{\alpha}$ denotes the Caputo
fractional derivative of order~$\alpha$ with respect to time $t$,
\begin{equation}\label{eq: Caputo}
(\fd{\alpha}v)(t)\defeq(\mathcal{I}^{1-\alpha}v')(t)
    =\int_0^t\omega_{1-\alpha}(t-s)v'(s)\zd{s}
    \quad\text{for $0<\alpha<1$ and $t>0$,}
\end{equation}
involving the Riemann--Liouville fractional integral operator of
order~$\beta>0$, defined by
\[
(\mathcal{I}^\beta v)(t)\defeq\int_0^t\omega_\beta(t-s)v(s)\zd{s}
    \quad\text{for $t>0$,}
    \quad\text{where $\omega_\beta(t)\defeq\frac{t^{\beta-1}}{\Gamma(\beta)}$.}
\]

An important and key consideration
\cite{FordYan,JinLazarovZhou:2016,JLZ17,McLeanMustapha:2007,
MustaphaAlMutawa:2012,MustaphaMcLean:2009,MustaphaMcLean:2013,StynesRiordanGracia:2016b}
in solving subdiffusion problems is that the solution~$u(x,t)$ is
typically non-smooth near the initial time, i.e., $\partial
u/\partial t=\mathcal{O}(1+t^{\alpha-1})$ as~$t\to0$, see \cite{McLean2010,SakamotoYamamoto2011,Stynes2016}. Among many
aproaches, one way to handle initial time singularity is to use
nonuniform time steps, see
\cite{Brunner1985,Brunner:book,McLean2010,McLeanMustapha:2007,LiaoZhaoTeng:2016,
LiaoLiZhang:2018,LiaoMcLeanZhang:2018,LiaoMcLeanZhang:2018b,LiaoYanZhang:2018,
MustaphaAlMutawa:2012,MustaphaMcLean:2009,MustaphaMcLean:2013,StynesRiordanGracia:2016b,ZhangSunLiao:2014}.
The main reason is that the nonuniform mesh is simple and flexible
to deal with not only the singular behavior near the initial time,
but also the possible rapid growth of the solution far away from $t
= 0$.

For the classical parabolic equation, the numerical analysis of the
widespread backward Euler and Crank-Nicolson schemes on general
nonuniform meshes for approximating the first-order time derivative
would be almost the same as the uniform case, and has been well
understand. For the subdiffusion problems considered here, the
numerical analysis on
nonuniform meshes is much complicate due to the convolution
integral form of Caputo  derivative \eqref{eq: Caputo}. Recently, Liao et al. developed a theoretical
framework
in~\cite{LiaoLiZhang:2018,LiaoMcLeanZhang:2018,LiaoMcLeanZhang:2018b,LiaoYanZhang:2018}
for the numerical analysis of nonuniform time approximations,
including the L1 formula
\cite{JinLazarovZhou:2016,LinXu:2007,LiaoLiZhang:2018,SunWu:2006},
two-level fast L1 formula \cite{LiaoYanZhang:2018} and the
fractional Crank-Nicolson (FracCN) scheme
\cite{Alikhanov2015,LiaoZhaoTeng:2016,LiaoMcLeanZhang:2018b},
to reaction-subdiffusion problems. This framework involves
three novel tools: a complementary discrete convolution kernel, a
discrete fractional Gr\"{o}nwall inequality and a global consistency
analysis. The stability and sharp $L^2$-norm error estimates are
obtained on general nonuniform meshes by taking into  the
initial singularity account. However, it seems that the framework is not
straightfoward to obtain the optimal $H^1$-norm estimates of
nonuniform time discretizations for problem~\eqref{eq: IBVP}.  This
motivates us to extend the framework to deal with the optimal
$H^1$-norm error estimate in this paper.

Actually, due to the nonlocal property of fractional time derivative
and the lack of smoothness near the initial time, the traditional
$H^1$-norm analysis (for the parabolic problems corresponding to
$\alpha\to1$) always leads to a suboptimal $H^1$-norm error
estimate. The goal of this paper is to achieve the optimal
$H^1$-norm error estimates of both L1 and FracCN schemes on a general
nonuniform mesh $0=t_0<t_1<t_2<\cdots<t_N=T$. Denote the time-step
size~$\tau_k\defeq t_k-t_{k-1}$, the adjoint step ratio
$\rho_k\defeq\tau_k/\tau_{k+1}$ for $k\ge1$, and the maximum
step size $\tau\defeq\max_{1\leq k\leq N}\tau_k$. Our focus is
on the time discretization of problems \eqref{eq: IBVP},  for
simplicity,  we only consider the finite difference method for the
spatial discretization in one dimension with $\Omega:=(x_l,x_r)$,
$$\opL:=-\partial_x(\mu(x)\partial_x)\quad\text{and}\quad
\text{$0<\mu_0\leq \mu(x)\le \mu_1$ for two positive constants
$\mu_0\ \text{and}\ \mu_1$.}$$

Nevertheless, the theoretical results in time approximations
together with their proofs here are also valid for multi-dimensional
problems, and are extendable for some other spatial discretization
such as the spectral method. To make the present analysis extendable
(such as for multi-term subdiffusion equations in Caputo's sense), let
$\sigma\in(0,1)\cup(1,2)$ be a regularity parameter and assume that
$u(\cdot,t)\in \mathrm{C}_{\sigma}^{3}((0,T])$, where the space
$\mathrm{C}_{\sigma}^{m}((0,T])$ is defined by
\begin{align}\label{function space}
\mathrm{C}_{\sigma}^{m}((0,T])\defeq\Big\{u\,\Big|\,& u\in  \mathrm{C}([0,T])\cap \mathrm{C}^{m}((0,T])\;\;\text{and}\nonumber\\
&\,\absb{u^{(\ell)}}\leq C_u\brab{1+t^{\sigma-\ell}}\;\;\text{for $\ell=1,\cdots, m$ and $0<t\leq T$}\Big\}.
\end{align}
Generally, the convergence rates of numerical Caputo derivatives are always limited by
the non-smoothness near the initial time.
It is reasonable to use a nonuniform mesh that concentrates grid points near $t=0$.
Let $\gamma\ge1$ be a user-chosen parameter, and assume that \cite{LiaoMcLeanZhang:2018b,LiaoYanZhang:2018,McLeanMustapha:2007}
\begin{description}
\item[\Mss{-conv}.]  There is a constant $C_{\gamma}>0$, independent of $k$,  such that
$\tau_k\le C_{\gamma}\tau\min\{1,t_k^{1-1/\gamma}\}$
for~$1\le k\le N$, $t_{k}\leq C_{\gamma}t_{k-1}$ and $\tau_k/t_k\le C_{\gamma}\tau_{k-1}/t_{k-1}$ for~$2\le k\le N$.
\end{description}
Since $\tau_1=t_1$, \Mss{-conv} implies that
$\tau_1=\mathcal{O}(\tau^{\gamma})$, while for those $t_k$ bounded
away from $t=0$ one has $\tau_k=\mathcal{O}(\tau)$. The
parameter~$\gamma$ controls the extent to which the grid points are
concentrated near~$t=0$. A practical example satisfying \Mss{-conv}
is an initially graded grid
\cite{Brunner1985,Brunner:book,McLeanMustapha:2007,MustaphaAlMutawa:2012,StynesRiordanGracia:2016b,LiaoLiZhang:2018}
\begin{equation}\label{eq: graded mesh}
t_k=(k/N_0)^\gamma T_0\;\;\text{for $0\leq k\leq N_0$}\quad
\text{and}\quad t_k=T_0+(k-N_0)\tau\;\;\text{for $N_0<k\le N$,}
\end{equation}
with
$$N_0\defeq\left\lceil\frac{\gamma NT_0}{T+(\gamma-1)T_0}\right\rceil\text{ for a small user-chosen $T_0\leq T$.}$$

Throughout the paper, any subscripted $C$, such as $C_\Omega$,
$C_{\gamma}$, $C_v$ and $C_u$, denotes a generic positive constant,
not necessarily the same at different occurrences, which is always
dependent on the given data and the solution, but independent of
temporal and spatial mesh sizes. The rest of this paper is organized
as follows. In Section~\ref{sec:general approach}, we present an
unified implicit time-stepping approach for subdiffusion problems
and some preliminary results. In Section~\ref{sec:L1 approach}, we
investigate the $H^1$-norm error bound for the L1 scheme, while the
second-order FracCN scheme with unequal time-steps is studied in
Section~\ref{sec:Alikhanov approach}. Two numerical examples in
Section~\ref{sec:numerical} are given to demonstrate the sharpness
of our analysis.

\section{An unified time-stepping scheme and $H^1$-norm stability}\label{sec:general approach}
\setcounter{equation}{0}

Assume that approximate the Laplacian~$\opL$ by the usual
second-order difference operator~$\opL_h$ on a discrete grid
$\overline{\Omega}_h\defeq\{\,x_l+ih\,|\,0\le i\le M\,\}$ with
$h\defeq (x_r-x_l)/M$. For any function $v_h$ on
$\overline{\Omega}_h$, we define
$\partial_hv_h(x_{i-\frac12}):=\brab{v_h(x_i)-v_h(x_{i-1})}/h$ and
$$\brab{\opL_hv_h}(x_i):=-\partial_h\brab{\mu(x_i)\partial_hv_h(x_i)}=-\brab{\mu(x_{i+\frac12})\partial_hv_h(x_{i+\frac12})
-\mu(x_{i-\frac12})\partial_hv_h(x_{i-\frac12})}/h.$$ We put
$\Omega_h\defeq\overline{\Omega}_h\cap\Omega$~and
$\partial\Omega_h\defeq\overline{\Omega}_h\cap\partial\Omega$. For
any functions $v_h$~and $w_h$ belonging to the space $\mathbb{V}_h$
of grid functions that vanish on the boundary~$\partial\Omega_h$, we
introduce the discrete inner product $\iprod{v_h,w_h}\defeq
h\sum_{x\in\Omega_h}v_h(x)w_h(x)$, the $L_2$ norm
$\|v_h\|\defeq\sqrt{\iprod{v_h,v_h}}\,$ and the $H^1$ semi-norm
$$|v_h|_1\defeq\sqrt{\iprod{v_h,\opL_hv_h}}=\sqrt{h\sum_{x\in\Omega_h}\mu(x)\bra{\partial_hv_h(x)}^2}\,.$$
There exists a positive constant $C_{\Omega}$ only dependent on the
domain $\Omega$, the constants $\mu_0$ and $\mu_1$ such that
$\|v_h\|\leq C_{\Omega}|v_h|_1$. The $H^1$ semi-norm $|v_h|_1$ is
equivalent to the discrete $H^1$ norm
$\|v_h\|_1:=\sqrt{\|v_h\|^2+|v_h|_1^2}$. So, in general, the
estimates of $|v_h|_1$ are called $H^1$-norm estimates.

\subsection{A time-weighted difference scheme}

Let $\nu\in[0,1/2)$ be an offset parameter and denote
$t_{n-\nu}\defeq\nu t_{n-1}+(1-\nu)t_n$.
For any mesh function $v^k\approx v(t_k)$,
define $v^{k-\nu}\defeq\nu v^{k-1}+(1-\nu)v^k$ and
$\diff v^k\defeq v^k-v^{k-1}$ for $k\ge1$.
The Caputo  derivative \eqref{eq: Caputo} of the function~$v$ can always
be approximated by a convolution-like summation,
\begin{equation}\label{eq: discrete Caputo}
(\fd{\alpha} v)(t_{n-\nu})\approx(\dfd{\alpha}v)^{n-\nu}\defeq\sum_{k=1}^n
    A^{(n,\nu)}_{n-k}\diff v^k\quad\text{for $1\le n\le N$,}
\end{equation}
with the local consistence error
\begin{equation}\label{eq: Truncation discrete Caputo}
\Upsilon^{n-\nu}[v]:=(\fd{\alpha}v)(t_{n-\nu})-(\dfd{\alpha}v)^{n-\nu}
\quad\text{for $1\le n\le N$.}
\end{equation}
Here, the corresponding discrete kernels, writing as~$A^{(n,\nu)}_{n-k}$ to reflect the convolution structure
of the integral
in~\eqref{eq: Caputo}, will be determined later. 
Our discrete solution, $u^n_h\approx u(x,t_n)$
for~$x\in\overline{\Omega}_h$, is defined by a time-weighted time-stepping scheme
\begin{equation}\label{eq: discrete IBVP}
\begin{aligned}
(\dfd{\alpha}u_h)^{n-\nu}+\opL_h u_h^{n-\nu}
    &=c(x) u_h^{n-\nu}+f(x,t_{n-\nu})&
    &\text{for $x\in\Omega_h$ and $1\le n\le N$,}\\
    u^n_h&=u_b(x,t_n)&
    &\text{for $x\in\partial\Omega_h$ and $1\le n\le N$,}\\
    u^0_h&=u_0(x)&&\text{for $x\in\Omega_h$.}
\end{aligned}
\end{equation}
In this paper, we will focus on two different cases of
$\brab{\dfd{\alpha}v}^{n-\nu}$: one is the widespread L1 formula
\cite{LiaoLiZhang:2018,LinXu:2007,StynesRiordanGracia:2016b,SunWu:2006}
with $\nu=0$, and the other is the recently suggested nonuniform
Alikhanov formula \cite{LiaoMcLeanZhang:2018b} with
$\nu=\theta:=\alpha/2$. For simplicity, the above  scheme \eqref{eq:
discrete IBVP} is called the L1 and FracCN method, respectively,
corresponding to the offset parameter $\nu=0$ and $\nu=\theta$.

The present approach would be fit for general nonuniform time meshes and
applicable for any discrete fractional derivatives having the
form~\eqref{eq: discrete Caputo} provided $A^{(n,\nu)}_{n-k}$~satisfy three criteria:
\begin{description}
\item[A1.] The discrete kernels are monotone, that is, $A^{(n,\nu)}_{k-2}\ge A^{(n,\nu)}_{k-1}>0$ for~$2\leq k\leq n\leq N$,
and the first one is properly large so that $(1-2\nu)A^{(n,\nu)}_{0}-(1-\nu)A^{(n,\nu)}_{1}\ge0$ for $\nu\in[0,1/2)$.
\item[A2.] There is a constant~$\pi_A>0$,
$A^{(n,\nu)}_{n-k}\ge\frac{1}{\pi_A\tau_k}\int_{t_{k-1}}^{t_k}\omega_{1-\alpha}(t_n-s)\zd{s}$
for $1\le k\le n\le N$.
\item[A3.] There is a constant~$\rho>0$ such that the local step ratio $\rho_k\le\rho$ for~$1\le k\le N-1$.
\end{description}
As noted in \cite{LiaoMcLeanZhang:2018}, the assumptions \Ass{1}--\Ass{2}
on the discrete convolution kernels~$A^{(n,\nu)}_{n-k}$ are valid for the most
frequently used discrete Caputo derivatives, at least if assumption~\Ass{3} is
satisfied for appropriate~$\rho$. Actually, the local mesh parameter~$\rho$ in~\Ass{3} will also
appear in our discrete fractional Gr\"{o}nwall inequality and
the $H^1$-norm stability estimate.

\subsection{An $H^1$-norm stability}

Always, the $H^1$-norm stability and convergence analysis on (general) nonuniform meshes makes use of
a discrete fractional Gr\"{o}nwall inequality and
a global consistency analysis, which involve a complementary discrete convolution
kernel~$P^{(n,\nu)}_{n-k}$ introduced by Liao et al.~\cite{LiaoLiZhang:2018,LiaoMcLeanZhang:2018} and
having the identical property
\begin{equation}\label{eq: P A}
\sum_{j=k}^nP^{(n,\nu)}_{n-j}A^{(j,\nu)}_{j-k}\equiv1\quad\text{for $1\le k\le n\le N$.}
\end{equation}
In fact, rearranging this identity yields a recursive formula (in
effect, a definition)
\begin{align}\label{eq: discreteConvolutionKernel-RL}
P_{0}^{(n,\nu)}\defeq\frac{1}{A_0^{(n,\nu)}},\quad
P_{n-j}^{(n,\nu)}\defeq
\frac{1}{A_0^{(j,\nu)}}
\sum_{k=j+1}^{n}\braB{A_{k-j-1}^{(k,\nu)}-A_{k-j}^{(k,\nu)}}P_{n-k}^{(n,\nu)}
    \quad \text{for $1\leq j\leq n-1$.}
\end{align}
Actually, it has been shown \cite[Lemma~2.2]{LiaoMcLeanZhang:2018} that $P^{(n,\nu)}_{n-k}$ is
well-defined and non-negative if the assumption \Ass{1} holds.
Furthermore, if the assumption \Ass{2} holds, then
\begin{align}\label{eq: P bound}
\sum_{j=1}^nP^{(n,\nu)}_{n-j}\omega_{1+m\alpha-\alpha}(t_n)\leq \pi_A\omega_{1+m\alpha}(t_n)\quad\text{for $1\le n\le N$ and $m=0,1$.}
\end{align}
Next we give a discrete fractional Gr\"{o}nwall inequality
\cite[Theorem 3.4]{LiaoMcLeanZhang:2018}, which should be fit for
the classical $H^1$-norm stability analysis.
\begin{theorem}\label{thm: H1 gronwall}
Let the assumptions \Ass{1}--\Ass{3} hold, let $0\le\nu<1/2$, and
let $(g^n)_{n=1}^N$ and $(\lambda_l)_{l=0}^{N-1}$ be given
non-negative sequences.  Assume further that there exists a constant
$\Lambda$ (independent of the time-step sizes) such that
$\Lambda\ge\sum_{l=0}^{N-1}\lambda_l$, and that the maximum step
size satisfies
\[
\taumax
    \le\frac{1}{\sqrt[\alpha]{2\pi_A\Gamma(2-\alpha)\Lambda}}\,.
\]
For any non-negative sequence~$(v^k)_{k=0}^N$ such that
\begin{equation}\label{eq: H1 Gronwall condition}
\sum_{k=1}^nA^{(n,\nu)}_{n-k}\diff v^k\le
    \sum_{k=1}^n\lambda_{n-k} v^{k-\nu}+g^n
    \quad\text{for $1\le n\le N$,}
\end{equation}
or
\begin{equation}\label{eq: H1 Gronwall condition2}
 v^n\le v^0+
    \sum_{j=1}^nP^{(n,\nu)}_{n-j}\sum_{k=1}^j\lambda_{j-k} v^{k-\nu}+\sum_{j=1}^nP^{(n,\nu)}_{n-j}g^j
    \quad\text{for $1\le n\le N$,}
\end{equation}
then it holds that
\begin{align*}
v^n&\le2E_\alpha\brab{2\max\{1,\rho\}\pi_A \Lambda t_n^\alpha}
    \biggl(v^0+\max_{1\le k\le n}\sum_{j=1}^k P^{(k,\nu)}_{k-j}g^j\biggr)\nonumber\\
&\le2E_\alpha\brab{2\max\{1,\rho\}\pi_A \Lambda t_n^\alpha}
    \biggl(v^0+\pi_A\Gamma(1-\alpha)\max_{1\le k\le n}\big\{t_k^{\alpha}g^k\big\}\biggr)\quad\text{for $1\le n\le N$.}
\end{align*}
\end{theorem}

In the subsequent discrete energy approach, we also need the following lemma,
which can be verified by a similar proof of \cite[Lemma~4.1]{LiaoMcLeanZhang:2018}.
\begin{lemma}\label{lem: general Dv v}
If the condition \Ass{1} holds, the discrete Caputo formula \eqref{eq: discrete Caputo} satisfies
$$v^{n-\nu}\bra{\dfd{\alpha}v}^{n-\nu}
    \geq\frac12\sum_{k=1}^nA^ {(n,\nu)}_{n-k}\diff\brab{\abst{v^k}^2}\quad
\text{for~$1\le n\le N$.}$$
\end{lemma}

We now consider the stability of the unified scheme \eqref{eq:
discrete IBVP} by assuming that $u_b(x,t_n)=0$. By taking the inner
product of the first equation in~\eqref{eq: discrete IBVP}
with~$2\bra{\dfd{\alpha}u_h}^{n-\nu}$, one has
\begin{align*}
2\mynormb{\bra{\dfd{\alpha}u_h}^{n-\nu}}^2
+2\myinnerb{\opL_hu_h^{n-\nu},\bra{\dfd{\alpha}u_h}^{n-\nu}}
&=2\myinnerb{cu_h^{n-\nu},\bra{\dfd{\alpha}u_h}^{n-\nu}}+2\myinnerb{f^{n-\nu},\bra{\dfd{\alpha}u_h}^{n-\nu}}\\
&\leq2\mynormb{\bra{\dfd{\alpha}u_h}^{n-\nu}}^2+\kappa^2\mynormb{u_h^{n-\nu}}^2+\mynormb{f^{n-\nu}}^2\,.
\end{align*}
Therefore, applying Lemma \ref{lem: general Dv v} ($v:=\opL_h^{1/2}u_h$) and the embedding inequality, one gets
\begin{align*}
\sum_{k=1}^nA^{(n,\nu)}_{n-k}\diff\brab{\absb{u_h^k}_1^2}
    \le&2\kappa^2C_{\Omega}\braB{(1-\nu)^2\absb{u_h^n}_1^2+\nu^2\absb{u_h^{n-1}}_1^2}
    +\mynormb{f^{n-\nu}}^2\\
    \le&2\kappa^2C_{\Omega}\braB{(1-\nu)\absb{u_h^n}_1^2+\nu\absb{u_h^{n-1}}_1^2}
    +\mynormb{f^{n-\nu}}^2,\quad1\le n\le N,
\end{align*}
which has the form of \eqref{eq: H1 Gronwall condition} with $\lambda_l=0$ for $l\ge1$,
$$\lambda_0:=2\kappa^2 C_{\Omega},\quad
v^k:=\absb{u_h^k}_1^2\quad\text{and}\quad
g^n:=\mynormb{f^{n-\nu}}^2.$$ Theorem \ref{thm: H1 gronwall} says
that the weighted time-stepping method~\eqref{eq: discrete IBVP} is
stable in the following sense.

\begin{theorem}\label{thm: stability}
If \Ass{1}--\Ass{3} hold with the maximum time-step size
$\taumax\le1/\sqrt[\alpha]{4\pi_A\Gamma(2-\alpha)\kappa^2C_{\Omega}}$,
then the time-stepping scheme~\eqref{eq: discrete IBVP}
with $u_b(x,t_n)=0$ is stable in the $H^1$-norm, that is,
\begin{align*}
\absb{u_h^n}_1^2&\le 2E_\alpha\bigl(4\pi_A\max\{1,\rho\}\kappa^2C_{\Omega}t_n^\alpha\bigr)\braB{\absb{u_h^{0}}_1^2
+\max_{1\le k\le n}\sum_{j=1}^k P^{(k,\nu)}_{k-j}\mynormb{f^{j-\nu}}^2\big\}}\\
&\le 2E_\alpha\bigl(4\pi_A\max\{1,\rho\}\kappa^2C_{\Omega}t_n^\alpha\bigr)\braB{\absb{u_h^{0}}_1^2
+\pi_A\Gamma(1-\alpha)\max_{1\le k\le n}\big\{t_k^{\alpha}\mynormb{f^{k-\nu}}^2\big\}}\;
\text{for $1\le n\le N.$}
\end{align*}
\end{theorem}

\subsection{An improved Gr\"{o}nwall inequality}

It is easy to check that, the solution error, $\tilde{u}^n_h:=u(x,t_n)-u^n_h$
for~$x\in\overline{\Omega}_h$, satisfies the zero-valued initial and boundary conditions, and the governing equation
\begin{equation}\label{eq: error discrete IBVP}
(\dfd{\alpha}\tilde{u}_h)^{n-\nu}+\opL_h\tilde{u}_h^{n-\nu}
    =c(x)\tilde{u}_h^{n-\nu}+\Upsilon_h^{n-\nu}[u]+R_w^{n-\nu}+R_s^{n-\nu}\quad
    \text{for $x\in\Omega_h$ and $1\le n\le N$,}
\end{equation}
where $\Upsilon_h^{n-\nu}[u]$ is defined by \eqref{eq: Truncation discrete Caputo},
\begin{equation}\label{eq: space truncation}
R^{n-\nu}_w\defeq\brab{c(x)-\opL}\kbrab{u^{n-\nu}-u(t_{n-\nu})}\quad\text{and}\quad
R^{n-\nu}_s\defeq(\opL-\opL_h)u^{n-\nu}\,.
\end{equation}
Nonetheless, it always yields a suboptimal $H^1$-norm error estimate if the a priori estimate in Theorem \ref{thm: stability}
is directly applied to the above error system, because the global consistency error
$\sum_{j=1}^k P^{(k,\nu)}_{k-j}\mynormb{\Upsilon_h^{j-\nu}[u]}^2$ has a loss of time accuracy, see an example in the next section.

To end this section, we present an extension of the fractional
Gr\"{o}nwall inequality in \cite[Theorem 3.1]{LiaoMcLeanZhang:2018}.
This result will be useful to obtain the optimal time accuracy in
the $H^1$-norm.
\begin{theorem}\label{thm: improved gronwall}
Let the assumptions \Ass{1}--\Ass{3} hold, let $0\le\nu<1/2$, and
let $(\xi^n)_{n=1}^N$, $(\eta^n)_{n=1}^N$ and
$(\lambda_l)_{l=0}^{N-1}$ be given non-negative sequences.  Assume
further that there exists a constant $\Lambda$ (independent of the
step sizes) such that $\Lambda\ge\sum_{l=0}^{N-1}\lambda_l$, and
that the maximum step size satisfies
\[
\taumax
    \le\frac{1}{\sqrt[\alpha]{2\pi_A\Gamma(2-\alpha)\Lambda}}\,.
\]
For any non-negative sequence~$(v^k)_{k=0}^N$ such that
\begin{equation}\label{eq: first Gronwall}
\sum_{k=1}^nA^{(n,\nu)}_{n-k}\diff\brab{v^k}^2\le
    \sum_{k=1}^n\lambda_{n-k}\brab{v^{k-\nu}}^2+v^{n-\nu}\xi^n+\bra{\eta^n}^2
    \quad\text{for $1\le n\le N$,}
\end{equation}
or
\begin{equation}\label{eq: transform first Gronwall}
\brab{v^n}^2\le\brab{v^0}^2
    +\sum_{j=1}^nP^{(n,\nu)}_{n-j}\sum_{k=1}^j\lambda_{j-k}\brab{v^{k-\nu}}^2+
\sum_{j=1}^nP^{(n,\nu)}_{n-j}v^{j-\nu}\xi^j+\sum_{j=1}^nP^{(n,\nu)}_{n-j}\bra{\eta^j}^2.
\end{equation}
Then it holds that,  for $1\le n\le N$,
\begin{equation}\label{eq: gronwall conclusion}
v^n\le2E_\alpha\brab{2\max\{1,\rho\}\pi_A \Lambda t_n^\alpha}
    \biggl(v^0+\max_{1\le k\le n}\sum_{j=1}^k P^{(k,\nu)}_{k-j}\xi^j+\sqrt{\pi_A\Gamma(1-\alpha)}\max_{1\le k\le n}\big\{t_k^{\alpha/2}\eta^k\big\}
    \biggr).
\end{equation}
\end{theorem}
\begin{proof}Two different cases are considered with a notation $\mathcal{E}_\alpha^n:=2E_\alpha\brab{2\max\{1,\rho\}\pi_A \Lambda t_n^\alpha}$.
If $$v^n\leq \eta^{*}:=\sqrt{\pi_A\Gamma(1-\alpha)}\max_{1\le k\le n}\big\{t_k^{\alpha/2}\eta^k\big\},$$
then the claimed inequality \eqref{eq: gronwall conclusion} follows because $\mathcal{E}_\alpha^n\ge2$ for any $0<\alpha<1$ and $n\ge0$. Otherwise,
if $v^n>\eta^{*}$,
then $v^n>\sqrt{\pi_A\Gamma(1-\alpha)}\,t_n^{\frac{\alpha}2}\eta^n$ and the inequality \eqref{eq: first Gronwall} becomes
\begin{equation}
\sum_{k=1}^nA^{(n,\nu)}_{n-k}\diff\brab{v^k}^2\le
    \sum_{k=1}^n\lambda_{n-k}\brab{v^{k-\nu}}^2+v^{n-\nu}\xi^n+v^n\frac{\eta^n}{\sqrt{\pi_A\Gamma(1-\alpha)t_n^{\alpha}}}
    \quad\text{for $1\le n\le N$.}
\end{equation}
Therefore, following the proof of \cite[Theorem 3.1]{LiaoMcLeanZhang:2018} with
$$g^n=\xi^n+\frac{\eta^n}{\sqrt{\pi_A\Gamma(1-\alpha)t_n^{\alpha}}},$$
one can apply \eqref{eq: P bound} to obtain that
\begin{align*}
v^n&\le\mathcal{E}_\alpha^n
    \biggl(v^0+\max_{1\le k\le n}\sum_{j=1}^k P^{(k,\nu)}_{k-j}\xi^j+\max_{1\le k\le n}\sum_{j=1}^k P^{(k,\nu)}_{k-j}\frac{\eta^j}{\sqrt{\pi_A\Gamma(1-\alpha)t_j^{\alpha}}}\biggr)\\
&\le\mathcal{E}_\alpha^n
    \biggl(v^0+\max_{1\le k\le n}\sum_{j=1}^k P^{(k,\nu)}_{k-j}\xi^j+\sqrt{\Gamma(1-\alpha)/\pi_A}\max_{1\le k\le n}\big\{t_k^{\alpha/2}\eta^k\big\}\max_{1\le k\le n}\sum_{j=1}^k P^{(k,\nu)}_{k-j}\omega_{1-\alpha}(t_j)\biggr)\\
&\le\mathcal{E}_\alpha^n
    \biggl(v^0+\max_{1\le k\le n}\sum_{j=1}^k P^{(k,\nu)}_{k-j}\xi^j
+\sqrt{\pi_A\Gamma(1-\alpha)}\max_{1\le k\le n}\big\{t_k^{\alpha/2}\eta^k\big\}\biggr).
\end{align*}
It completes the proof.
\end{proof}

\begin{remark}\label{rem: sum P}
One may use the inequality \eqref{eq: P bound} to bound the
summation $\displaystyle\sum_{j=1}^k P^{(k,\nu)}_{k-j}\xi^j$, that
is,
\begin{equation*}
\sum_{j=1}^k P^{(k,\nu)}_{k-j}\xi^j\leq \sum_{j=1}^k
P^{(k,\nu)}_{k-j}\omega_{1-\alpha}(t_j)\max_{1\leq j\leq
k}\frac{\xi^j}{\omega_{1-\alpha}(t_j)} \leq \pi_A\max_{1\leq j\leq
k}\frac{\xi^j}{\omega_{1-\alpha}(t_j)}.
\end{equation*}
So the discrete solution of \eqref{eq: first Gronwall} can also be bounded by
\[
v^n\le2E_\alpha\brab{2\max\{1,\rho\}\pi_A \Lambda t_n^\alpha}
    \braB{v^0+\pi_A\Gamma(1-\alpha)\max_{1\le j\le n}\{t_j^{\alpha}\xi^j\}+
\sqrt{\pi_A\Gamma(1-\alpha)}\max_{1\le k\le n}\big\{t_k^{\alpha/2}\eta^k\big\}}.
\]
On the other hand, if the given sequence $(\lambda_l)_{l=0}^{N-1}$ is non-positive and the constant $\Lambda\le0$,
a similar argument will show that the discrete inequality \eqref{eq: gronwall conclusion} holds in a simpler form,
requiring only the assumptions \Ass{1}-\Ass{2} but no restrictions on time steps,
\begin{align}\label{eq: gronwall conclusion simple}
v^n&\le v^0+\max_{1\le k\le n}\sum_{j=1}^k P^{(k,\nu)}_{k-j}\xi^j
+\sqrt{\pi_A\Gamma(1-\alpha)}\max_{1\le k\le n}\big\{t_k^{\alpha/2}\eta^k\big\}\nonumber\\
&\le v^0+\pi_A\Gamma(1-\alpha)\max_{1\le j\le n}\{t_j^{\alpha}\xi^j\}
+\sqrt{\pi_A\Gamma(1-\alpha)}\max_{1\le k\le n}\big\{t_k^{\alpha/2}\eta^k\big\}
    \quad\text{for $1\le n\le N$.}
\end{align}
\end{remark}

\section{Sharp $H^1$-norm error estimate for L1 scheme}\label{sec:L1 approach}
\setcounter{equation}{0}

In this section, assume that the solution $u\in \mathrm{C}([0,T];H^4(\Omega))\cap\mathrm{C}_{\sigma}^{2}((0,T];H^1(\Omega))$.
The Caputo's derivative $\fd{\alpha}v$ is approximated by the L1 formula
$(\dfd{\alpha}v)^n$, the case of $\nu=0$ in \eqref{eq: discrete Caputo}, with unequal time-steps.
The corresponding discrete convolution kernel $A_{n-k}^{(n,0)}$ is defined by
\begin{align}\label{Coefficient1-L1formula}
A_{n-k}^{(n,0)}\defeq\int_{t_{k-1}}^{t_{k}}\frac{\omega_{1-\alpha}(t_n-s)}{\tau_k}\zd s
=\frac{\omega_{2-\alpha}(t_n-t_{k-1})-\omega_{2-\alpha}(t_n-t_{k})}{\tau_k}\,,\quad 1\leq k\leq n.
\end{align}
Obviously, \Ass{2} holds for $\pi_A=1$, and next Lemma implies that
\Ass{1} is valid.

\begin{lemma}\cite[Lemma 2.1]{LiaoYanZhang:2018}\label{lemma: L1 kernel}
For fixed $n\geq2$, the discrete kernel $A_{n-k}^{(n,0)}$ in
\eqref{Coefficient1-L1formula} satisfies
$$A^{(n,0)}_{n-k-1}-A^{(n,0)}_{n-k}>\frac12\int_{t_{k-1}}^{t_k}\zd\omega_{1-\alpha}(t_{n}-s)>0\quad\text{for $1\leq k\leq n-1$.}$$
\end{lemma}

Hence we can use the complementary discrete convolution
kernel~$P^{(n,0)}_{n-k}$, see \eqref{eq: P A}-\eqref{eq: P bound}, in the subsequent analysis.
Also, Lemma \ref{lemma: L1 kernel} and Theorem \ref{thm: stability} imply
the unconditional stability of L1 scheme for the linear problem \eqref{eq: IBVP}.
\begin{corollary}\label{corol: L1 stability}
The L1 method \eqref{eq: discrete IBVP} with $\nu=0$ is stable in the discrete $H^1$ norm.
\end{corollary}

For the L1 scheme \eqref{eq: discrete Caputo} with the discrete convolution kernels \eqref{Coefficient1-L1formula},
we have the following estimate on the consistency error.
\begin{lemma}\label{lemma:L1formulaNonuniform-consistence}
For $v\in \mathrm{C}_{\sigma}^{2}((0,T])$, the local consistency
error of the L1 formula $(\dfd{\alpha}v)^n$ satisfies
\begin{align*}
\absb{\Upsilon^{n}[v]}\leq
A_{0}^{(n,0)}G^n+\sum_{k=1}^{n-1}\brab{A_{n-k-1}^{(n,0)}-A_{n-k}^{(n,0)}}G^k\quad\text{for~$n\ge 1$,}
\end{align*}
where $G^k$ is defined by
$G^k:=2\displaystyle\int_{t_{k-1}}^{t_k}\bra{t-t_{k-1}}\abs{v''(t)}\zd
t.$ Thus the global consistency error
\begin{align*}
\sum_{j=1}^nP^{(n,0)}_{n-j}\absb{\Upsilon^{j}[v]}\leq
\frac{C_v}{\sigma}\tau_1^{\sigma}+
\frac{C_v}{1-\alpha}\max_{2\leq k\leq n}t_{k}^{\alpha}t_{k-1}^{\sigma-2}\tau_k^{2-\alpha}\quad\text{for~$n\ge 1$.}
\end{align*}
Moreover, if the time mesh satisfies \Mss{-conv}, then
\begin{align*}
\sum_{j=1}^nP^{(n,0)}_{n-j}\absb{\Upsilon^{j}[v]}\leq
\frac{C_v}{\sigma(1-\alpha)}\tau^{\min\{2-\alpha,\gamma\sigma\}}\quad\text{for~$1\le n\le N$.}
\end{align*}
\end{lemma}
\begin{proof}
 See the proof of Lemmas 3.1 and 3.3 (taking
$\epsilon=0$) in \cite{LiaoYanZhang:2018}.
\end{proof}

\subsection{Suboptimal estimate by traditional $H^1$-norm analysis}
In this subsection, we show that the traditional $H^1$-norm analysis
together with the discrete Gr\"{o}nwall inequality in Theorem
\ref{thm: H1 gronwall} always yields a suboptimal estimate in the
$H^1$-norm, if the solution is nonsmooth near the initial time.
Without losing the generality, we consider the error equation
\eqref{eq: error discrete IBVP} with $\nu=0$, that is,
\begin{equation}\label{eq: simple L1 error eqn}
(\dfd{\alpha}\tilde{u}_h)^{n}+\opL_h\tilde{u}_h^{n}
    =c(x)\tilde{u}_h^{n}+\Upsilon_h^{n}[u]+R_s^{n}\quad
    \text{for $x\in\Omega_h$ and $1\le n\le N$,}
\end{equation}
where $\Upsilon_h^{n}[u]$ and $R^{n}_s$ are defined by \eqref{eq: Truncation discrete Caputo} and \eqref{eq: space truncation}, respectively.
Taking the inner product of the error equation
in~\eqref{eq: simple L1 error eqn} with~$2(\dfd{\alpha}\tilde{u}_h)^{n}$, one has
\begin{align*}
2\mynormb{\bra{\dfd{\alpha}\tilde{u}_h}^{n}}^2+2\myinnerb{\opL_h\tilde{u}_h^{n},(\dfd{\alpha}\tilde{u}_h)^{n}}
=2\myinnerb{c \tilde{u}^{n}_h,(\dfd{\alpha}\tilde{u}_h)^{n}}
+2\myinnerb{\Upsilon_h^n[u]+R_s^n,(\dfd{\alpha}\tilde{u}_h)^{n}}\\
\leq
2\mynormb{\bra{\dfd{\alpha}\tilde{u}_h}^{n}}^2+\kappa^2\mynormb{\tilde{u}_h^{n}}^2+
\mynormb{\Upsilon_h^n[u]+R_s^n}^2\,.
\end{align*}
Lemma \ref{lemma: L1 kernel} ensures \Ass{1}, so we apply Lemma \ref{lem: general Dv v} and the embedding inequality to get
\begin{align*}
\sum_{k=1}^nA^{(n,0)}_{n-k}\diff\brab{\absb{\tilde{u}_h^{n}}_1^2}
    \le&\,\kappa^2C_{\Omega}\absb{\tilde{u}_h^{n}}_1^2
    +\mynormb{\Upsilon_h^n[u]+R_s^n}^2,
\end{align*}
which takes the form of \eqref{eq: H1 Gronwall condition} with
$v^k=\absb{\tilde{u}_h^{k}}_1^2$ and
$g^n=\mynormb{\Upsilon_h^n[u]+R_s^n}^2$. So Theorem \ref{thm: H1
gronwall} together with the upper bound \eqref{eq: P bound} yields
the following estimate
\begin{align}\label{estimate: L1 standard H1 norm}
\absb{\tilde{u}_h^{n}}_1^2\le
2E_\alpha\brab{2\max(1,\rho)\kappa^2C_{\Omega}t_n^\alpha}
    \biggl(\max_{1\le k\le n}\sum_{j=1}^k P^{(k,0)}_{k-j}\mynormb{\Upsilon_h^j[u]}^2
+\Gamma(1-\alpha)\max_{1\le k\le n}t_k^{\alpha}\mynormb{R_s^k}^2
\biggr),
\end{align}
if the maximum time-step size
$\taumax\le1/\sqrt[\alpha]{2\Gamma(2-\alpha)\kappa^2C_{\Omega}}$. To
continue the error analysis, one requires the following result,
which takes advantage of the discrete convolution structure of local
truncation error in Lemma
\ref{lemma:L1formulaNonuniform-consistence}.


\begin{lemma}\label{lemma:L1formulaNonuniform-H1consistence}
If $v\in \mathrm{C}_{\sigma}^{2}((0,T])$ for
$\sigma\in(\frac{\alpha}2,1)\cup(1,2)$ and the maximum step ratio
$\rho\leq1$, then
\begin{align*}
\sum_{j=1}^nP^{(n,0)}_{n-j}\absb{\Upsilon^j[v]}^2\leq
\frac{C_v}{\sigma^2}\tau_1^{2\sigma-\alpha}
+\frac{C_v}{1-\alpha}\max_{2\leq k\leq n}t_{k}^{\alpha}t_{k-1}^{2\sigma-4}\tau_k^{4-2\alpha}
\quad\text{for~$1\le n\le N$.}
\end{align*}
Moreover, if the time mesh satisfies \Mss{-conv}, then
\begin{align*}
\sum_{j=1}^nP^{(n,0)}_{n-j}\absb{\Upsilon^j[v]}^2\leq
\frac{C_v}{\sigma^2(1-\alpha)}\tau^{2\min\{2-\alpha,\gamma(\sigma-\alpha/2)\}}\quad\text{for~$1\le n\le N$.}
\end{align*}
\end{lemma}
\begin{proof}
Applying Lemma \ref{lemma:L1formulaNonuniform-consistence} and the Cauchy-Schwarz inequality, one has
\begin{align*}
\absb{\Upsilon^j[v]}^2&\leq\brab{2A_{0}^{(j,0)}-A_{j-1}^{(j,0)}}
\kbraB{A_{0}^{(j,0)}(G^j)^2+\sum_{k=1}^{j-1}\brab{A_{j-k-1}^{(j,0)}-A_{j-k}^{(j,0)}}(G^k)^2}\\
&\leq
2\brab{A_{0}^{(j,0)}G^j}^2+2\sum_{k=1}^{j-1}A_{0}^{(j,0)}\brab{A_{j-k-1}^{(j,0)}-A_{j-k}^{(j,0)}}(G^k)^2
\end{align*}
The definition \eqref{Coefficient1-L1formula} gives
$A_{0}^{(j,0)}=\tau_j^{-\alpha}/\Gamma(2-\alpha)$ such that
$\displaystyle\max_{k+1\leq j\leq n}A_{0}^{(j,0)}=A_{0}^{(k,0)}$ if
the maximum ratio $\rho\leq1$. Multiplying the above inequality
by~$P^{(n,0)}_{n-j}$ and summing the index $j$ from $1$ to $n$, we
exchange the order of summation and apply the definition \eqref{eq:
discreteConvolutionKernel-RL} of $P_{n-j}^{(n,0)}$ to get
\begin{align}\label{globalApproximateError-immediate}
\sum_{j=1}^nP_{n-j}^{(n,0)}\absb{\Upsilon^j[v]}^2
\leq&\,
2\sum_{j=1}^nP_{n-j}^{(n,0)}\brab{A_{0}^{(j,0)}G^j}^2+
2\sum_{j=2}^nP_{n-j}^{(n,0)}\sum_{k=1}^{j-1}A_{0}^{(j,0)}\brab{A_{j-k-1}^{(j,0)}-A_{j-k}^{(j,0)}}(G^k)^2
\nonumber\\
=&\,
2\sum_{j=1}^nP_{n-j}^{(n,0)}\brab{A_{0}^{(j,0)}G^j}^2+
2\sum_{k=1}^{n-1}(G^k)^2\sum_{j=k+1}^nP_{n-j}^{(n,0)}A_{0}^{(j,0)}\brab{A_{j-k-1}^{(j,0)}-A_{j-k}^{(j,0)}}
\nonumber\\
\leq&\,
2\sum_{j=1}^nP_{n-j}^{(n,0)}\brab{A_{0}^{(j,0)}G^j}^2+
2\sum_{k=1}^{n-1}(G^k)^2A_{0}^{(k,0)}\sum_{j=k+1}^nP_{n-j}^{(n,0)}\brab{A_{j-k-1}^{(j,0)}-A_{j-k}^{(j,0)}}
\nonumber\\
=&\,
2\sum_{j=1}^nP_{n-j}^{(n,0)}\brab{A_{0}^{(j,0)}G^j}^2+
2\sum_{k=1}^{n-1}P_{n-k}^{(n,0)}\brab{A_{0}^{(k,0)}G^k}^2\quad\text{for~$1\le n\le N$.}
\end{align}
Now, following the proof of Lemma 3.3 in \cite{LiaoLiZhang:2018},
one can apply the definition \eqref{Coefficient1-L1formula} to find that
\begin{align*}
\frac{A_{0}^{(k,0)}}{A_{k-2}^{(k,0)}}<\frac{\omega_{2-\alpha}(\tau_k)}{\tau_k\,\omega_{1-\alpha}(t_k-t_1)}
=\frac{(t_k-t_1)^{\alpha}}{1-\alpha}\tau_k^{-\alpha},\quad 2\leq k\leq n\,.
\end{align*}
The regularity assumption implies that
$$G^1\leq C_v\tau_1^{\sigma}/\sigma\quad\text{and}\quad G^k\leq C_vt_{k-1}^{\sigma-2}\tau_k^2\quad \text{for $2\leq k\leq n.$}$$
Furthermore, the property \eqref{eq: P A} shows that
$P_{n-1}^{(n,0)}A_{0}^{(1,0)}\leq1$ and
$\displaystyle\sum_{k=2}^{n}P_{n-k}^{(n,0)}A_{k-2}^{(k,0)}=1.$ Thus
it follows from \eqref{globalApproximateError-immediate}  that
\begin{align*}
\sum_{j=1}^nP_{n-j}^{(n,0)}\absb{\Upsilon^j[v]}^2\leq&\,
4P_{n-1}^{(n,0)}\brab{A_{0}^{(1,0)}G^1}^2+4\sum_{k=2}^{n}P_{n-k}^{(n)}\brab{A_{0}^{(k,0)}G^k}^2\\
\leq&\,4A_{0}^{(1,0)}\brab{G^1}^2+\frac{4}{1-\alpha}\sum_{k=2}^{n}P_{n-k}^{(n,0)}A_{k-2}^{(k,0)}t_k^{\alpha}\tau_k^{-\alpha}A_{0}^{(k,0)}\brab{G^k}^2\\
\leq&\,
C_v\tau_1^{2\sigma-\alpha}/\sigma^2+\frac{C_v}{1-\alpha}\sum_{k=2}^{n}P_{n-k}^{(n)}A_{k-2}^{(k,0)}t_{k}^{\alpha}t_{k-1}^{2\sigma-4}\tau_k^{4-2\alpha}\\
\leq&\,
C_v\tau_1^{2\sigma-\alpha}/\sigma^2+\frac{C_v}{1-\alpha}\max_{2\leq k\leq n}t_{k}^{\alpha}t_{k-1}^{2\sigma-4}\tau_k^{4-2\alpha}\quad\text{for~$1\le n\le N$.}
\end{align*}
If the mesh fulfills \Mss{-conv}, then
$\tau_1\le C_{\gamma}\tau^{\gamma}$ and, with~$\beta:=2\min\{2-\alpha,\gamma(\sigma-\alpha/2)\}$,
\begin{align*}
t_{k}^{\alpha}t_{k-1}^{2\sigma-4}\tau_k^{4-2\alpha}&\leq C_{\gamma}t_{k}^{2\sigma-4+\alpha}\tau_k^{4-2\alpha}
\leq C_\gamma t_k^{2\sigma-4+\alpha}\tau_k^{4-2\alpha-\beta}
    \bigl(\tau\min\{1,t_k^{1-1/\gamma}\}\bigr)^{\beta}\\
&\leq C_\gamma t_k^{2\sigma-\alpha-\beta/\gamma}(\tau_k/t_k)^{4-2\alpha-\beta}\tau^{\beta}
\leq C_\gamma t_k^{\max\{0,2\sigma-\alpha-(4-2\alpha)/\gamma\}}\tau^{\beta}\quad\text{for~$2\le k\le n$.}
\end{align*}
It leads to the desired estimate and completes the proof.
\end{proof}

Since $u\in \mathrm{C}([0,T];H^4(\Omega))$, one has the spatial error estimate $\mynormb{R_s^k}\le C_uh^2$.
Applying the zero-valued initial data and Lemma \ref{lemma:L1formulaNonuniform-H1consistence},
one derive from \eqref{estimate: L1 standard H1 norm} that
\begin{align*}
\absb{\tilde{u}_h^{n}}_1^2&\le
C_u\max_{1\le k\le n}\sum_{j=1}^k P^{(k,0)}_{k-j}\mynormb{\Upsilon_h^j[u]}^2+C_u\Gamma(1-\alpha)\max_{1\le k\le n}t_k^{\alpha}\mynormb{R_s^k}^2\\
&\le
\frac{C_u}{\sigma^2(1-\alpha)}\brab{\tau^{2\min\{2-\alpha,\gamma(\sigma-\alpha/2)\}}+
t_n^{\alpha}h^4}
\end{align*}
or
\begin{align}\label{estimate: L1 suboptimal H1 norm}
\absb{\tilde{u}_h^{n}}_1\le
\frac{C_u}{\sigma\sqrt{1-\alpha}}\brab{\tau^{\min\{2-\alpha,\gamma(\sigma-\alpha/2)\}}+t_n^{\alpha/2}h^2}\quad\text{for $1\le n\le N$.}
\end{align}
It is optimal only when the regularity parameter
$\sigma\ge2-\alpha/2$, see previous studies
\cite{LiaoZhaoTeng:2016,ZhangSunLiao:2014} by assuming the solution
is smooth near the initial time; however, there is always a loss of
theoretical accuracy $\mathcal{O}(\tau^{-\gamma\alpha/2})$ in time
under the realistic assumption.

\begin{remark}\label{remark: L2 analysis}
As similar to the ordinary diffusion case corresponding to $\alpha=1$,
the standard $L^2$-norm error analysis \cite{LiaoLiZhang:2018} leads to the sharp estimate for
a weighted $H^1$ norm, but always gives a suboptimal estimate for the $H^1$-norm error at any time $t_n$,
see also \cite{LinXu:2007,LvXu:2016,SunWu:2006} for the analysis considering the smooth solutions.
Actually, by taking the inner product of the error equation
in~\eqref{eq: simple L1 error eqn} with~$\tilde{u}_h^{n}$ and applying the Cauchy-Schwarz inequality, one has
\begin{align*}
\myinnerb{(\dfd{\alpha}\tilde{u}_h)^{n},\tilde{u}_h^{n}}+\absb{\tilde{u}_h^{n}}_1^2
\leq
\kappa\mynormb{\tilde{u}^{n}_h}^2+\mynormb{\tilde{u}^{n}_h}\mynormb{\Upsilon_h^{n}[u]+R_s^n}\,.
\end{align*}
Therefore, applying Lemma \ref{lem: general Dv v}, we have
\begin{align*}
\sum_{k=1}^jA^{(j,0)}_{j-k}\diff\brab{\mynormb{\tilde{u}_h^k}^2}+2\absb{\tilde{u}_h^{j}}_1^2
\leq
2\kappa\mynormb{\tilde{u}^{j}_h}^2+2\mynormb{\tilde{u}^{j}_h}\mynormb{\Upsilon_h^{j}[u]+R_s^j}\,.
\end{align*}
Multiplying the above inequality by~$P^{(n,0)}_{n-j}$
and summing the index $j$ from $1$ to $n$, we get
\begin{align*}
\mynormb{\tilde{u}_h^n}^2+2\sum_{j=1}^nP^{(n,0)}_{n-j}\absb{\tilde{u}_h^{j}}_1^2
\leq\mynormb{\tilde{u}_h^0}^2+
2\kappa\sum_{j=1}^nP^{(n,0)}_{n-j}\mynormb{\tilde{u}^{j}_h}^2
+2\sum_{j=1}^nP^{(n,0)}_{n-j}\mynormb{\tilde{u}^{j}_h}\mynormb{\Upsilon_h^{j}[u]+R_s^j}\,,
\end{align*}
or, with $\absb{\!\absb{\!\absb{\tilde{u}_h^n}\!}\!}:=\sqrt{\mynormb{\tilde{u}_h^n}^2+2\sum_{j=1}^n P^{(n,0)}_{n-j}\absb{\tilde{u}_h^{j}}_1^2}$,
\begin{align*}
\absb{\!\absb{\!\absb{\tilde{u}_h^n}\!}\!}^2
    \le&\,\absb{\!\absb{\!\absb{\tilde{u}_h^0}\!}\!}^2+2\kappa \sum_{j=1}^n P^{(n,0)}_{n-j}\absb{\!\absb{\!\absb{\tilde{u}_h^j}\!}\!}^2
    +2\sum_{j=1}^n P^{(n,0)}_{n-j}\absb{\!\absb{\!\absb{\tilde{u}_h^j}\!}\!}\mynormb{\Upsilon_h^{j}[u]+R_s^j},
\end{align*}
which takes the form of \eqref{eq: transform first Gronwall} with
$v^n=\absb{\!\absb{\!\absb{\tilde{u}_h^n}\!}\!}$,
$\xi^j=2\mynormb{\Upsilon_h^{j}[u]+R_s^j}$ and $\eta^j=0$. So the
discrete Gr\"{o}nwall inequality in Theorem \ref{thm: improved
gronwall} yields
\begin{align*}
\absb{\!\absb{\!\absb{\tilde{u}_h^n}\!}\!}\le&\,4E_\alpha\brab{2\max(1,\rho)\kappa t_n^\alpha}
    \biggl(\mynormb{\tilde{u}_h^0}+\max_{1\le k\le n}\sum_{j=1}^k P^{(k,0)}_{k-j}\mynormb{\Upsilon_h^{j}[u]+R_s^j}\biggr),
\end{align*}
if the maximum step size
$\taumax\le1/\sqrt[\alpha]{4\Gamma(2-\alpha)\kappa}$. Applying the
spatial error estimate $\mynormb{R_s^k}\le C_uh^2$ and Lemma
\ref{lemma:L1formulaNonuniform-consistence}, we have the sharp
estimate for a weighted $H^1$-norm
\begin{align*}
\sqrt{\mynormb{\tilde{u}_h^n}^2+\sum_{j=1}^n P^{(n,0)}_{n-j}\absb{\tilde{u}_h^{j}}_1^2}
&\leq \frac{C_u}{\sigma(1-\alpha)}\brab{\tau^{\min\{2-\alpha,\gamma\sigma\}}+t_n^{\alpha}h^2}
\quad\text{for $1\le n\le N$.}
\end{align*}
However, a loss of accuracy $O(\tau_n^{-\frac{\alpha}2})$ will be seen in the $H^1$-norm error at any time $t_n$,
\begin{align*}
\absb{\tilde{u}_h^{n}}_1\leq \sqrt{A^{(n,0)}_{0}\sum_{j=1}^n P^{(n,0)}_{n-j}\absb{\tilde{u}_h^{j}}_1^2}
&\leq \frac{C_u}{\sigma(1-\alpha)}\tau_n^{-\frac{\alpha}2}
\brab{\tau^{\min\{2-\alpha,\gamma\sigma\}}+t_n^{\alpha}h^2}\quad\text{for $1\le n\le N$.}
\end{align*}
Compared with \eqref{estimate: L1 suboptimal H1 norm}, the loss of accuracy appears both in time and space.
\end{remark}

\subsection{Sharp $H^1$-norm error estimate}
A sharp $H^1$-norm error estimate reflecting the initial singularity is obtained by applying
the improved discrete Gr\"{o}nwall inequality in Theorem \ref{thm: improved gronwall}
and treating the temporal truncation error specially.
We will redefine the time truncation error uniformly over the closed space
domain, that is,
$\Upsilon_h^{n}[u]$ in \eqref{eq: simple L1 error eqn} can be redefined as follows, see Remark \ref{remark: H2 analysis} below,
\begin{align}
\overline{\Upsilon}_h^{n}[u]&\defeq (\fd{\alpha}u_h)(t_{n})-(\dfd{\alpha}u_h)^{n}
\quad\text{for $x\in\bar{\Omega}_h$ and $1\le n\le N$.}\label{eq: L1 time truncation}
\end{align}
Then the error equation \eqref{eq: simple L1 error eqn} can be formulated as
\begin{equation}\label{eq: L1 error eqn}
(\dfd{\alpha}\tilde{u}_h)^{n}=-\opL_h\tilde{u}_h^{n}
    +c(x)\tilde{u}_h^{n}+\overline{\Upsilon}_h^{n}[u]+R_s^n
    \quad\text{for $x\in\Omega_h$ and $1\le n\le N$.}
\end{equation}
By taking the inner product of the error equation
in~\eqref{eq: L1 error eqn} with~$\opL_h\tilde{u}_h^{n}$, one applies the discrete first Green formula to find
\begin{align*}
\myinnerb{\bra{\dfd{\alpha}\tilde{u}_h}^{n},\opL_h\tilde{u}_h^{n}}
=&\,-\mynormb{\opL_h\tilde{u}_h^{n}}^2
+\myinnerb{c\tilde{u}_h^{n},\opL_h\tilde{u}_h^{n}}
+\myinnerb{\mu\partial_h\overline{\Upsilon}_h^{n}[u],\partial_h\tilde{u}_h^{n}}
+\myinnerb{R_s^n,\opL_h\tilde{u}_h^{n}}\\
\leq&\,
\frac{1}2\kappa^2\mynormb{\tilde{u}_h^{n}}^2+
\absb{\tilde{u}_h^{n}}_1\absb{\overline{\Upsilon}_h^{n}[u]}_1+\frac1{2}\mynormb{R_s^n}^2\\
\leq&\,
\frac{1}2\kappa^2C_{\Omega}\absb{\tilde{u}_h^{n}}_1^2+
\absb{\tilde{u}_h^{n}}_1\absb{\overline{\Upsilon}_h^{n}[u]}_1+\frac1{2}\mynormb{R_s^n}^2,
\end{align*}
where the Cauchy-Schwarz inequality and the embedding inequality have been used.
We apply Lemma \ref{lem: general Dv v} to obtain
\begin{align*}
\sum_{k=1}^nA^{(n,0)}_{n-k}\diff\brab{\absb{\tilde{u}_h^{k}}_1^2}
    \le&\,\kappa^2C_{\Omega}\absb{\tilde{u}_h^{n}}_1^2
    +2\absb{\tilde{u}_h^{n}}_1\absb{\overline{\Upsilon}_h^{n}[u]}_1
    +\mynormb{R_s^n}^2,\quad1\le n\le N,
\end{align*}
which has the form of \eqref{eq: first Gronwall} with
$\lambda_0:=\kappa^2C_{\Omega}$, $v^k:=\absb{\tilde{u}_h^{k}}_1$,
$\xi^n:=2\absb{\overline{\Upsilon}_h^{n}[u]}_1$ and
$\eta^n:=\mynormb{R_s^n}$. Therefore, applying Theorem \ref{thm:
improved gronwall}, we see that
\begin{equation}\label{eq: L1 e R}
\absb{\tilde{u}_h^{n}}_1\le4E_\alpha\brab{2\max\{1,\rho\}\kappa^2C_{\Omega}t_n^\alpha}
    \max_{1\le k\le n}\biggl(\sum_{j=1}^k P^{(k,0)}_{k-j}\absb{\overline{\Upsilon}_h^{j}[u]}_1
  +\sqrt{\Gamma(1-\alpha)}t_k^{\alpha/2}\mynormb{R_s^k}
    \biggr),
\end{equation}
if the maximum time-step size
$\taumax\le1/\sqrt[\alpha]{2\Gamma(2-\alpha)\kappa^2C_{\Omega}}$. It
remains to evaluate the right-hand side of \eqref{eq: L1 e R} by
taking the initial singularity into account. Note that, the formula
of Taylor expansion with integral remainder gives
\begin{align*}
&\partial_h\brab{\overline{\Upsilon}_h^{n}[u]}(x_{i-\frac12})
=\int_0^1\brab{\fd{\alpha}\partial_xu(x_i-\lambda h)}(t_{n})
-\brab{\dfd{\alpha}\partial_xu(x_i-\lambda h)}^{n}\zd\lambda
\end{align*}
for $x\in\Omega_h$ and $1\le n\le N$. Then Lemma \ref{lemma:L1formulaNonuniform-consistence} with $v=\partial_xu$
gives the global consistency error
$$\sum_{j=1}^k P^{(k,0)}_{k-j}\absb{\overline{\Upsilon}_h^{j}[u]}_1\leq \frac{C_u}{\sigma(1-\alpha)}\tau^{\min\{2-\alpha,\gamma\sigma\}}\quad\text{for~$1\le k\le N$,}$$
and, obviously, $\mynormb{R_s^k}\leq C_uh^2$. So the inequality \eqref{eq: L1 e R} shows that
\begin{equation*}
\absb{\tilde{u}_h^{n}}_1\le \frac{C_u}{\sigma(1-\alpha)}E_\alpha\brab{2\max\{1,\rho\}\kappa^2C_{\Omega}t_n^\alpha}\bra{\tau^{\min\{2-\alpha,\gamma\sigma\}}
+t_n^{\alpha/2}h^2}\quad\text{for~$1\le n\le N$.}
\end{equation*}
It yields the following $H^1$-norm error estimate.
\begin{theorem}\label{th:Convergence-nonuniformL1Scheme}
Assume that the sub\/diffusion solution $u\in
\mathrm{C}([0,T];H^4(\Omega))\cap\mathrm{C}_{\sigma}^{2}((0,T];H^1(\Omega))$.
If the maximum time-step size
$\taumax\le1/\sqrt[\alpha]{2\Gamma(2-\alpha)\kappa^2C_{\Omega}}\,$,
then the solution of the L1 method \eqref{eq: discrete IBVP} with
$\nu=0$ on the nonuniform mesh satisfying \Ass{3} and \Mss{-conv},
is unconditionally convergent in the discrete $H^1$-norm,
\begin{align}\label{L1convegenceerror}
\absb{u(t_n)-u_h^n}_1\leq
\frac{C_u}{\sigma(1-\alpha)}E_\alpha\brab{2\max\{1,\rho\}\kappa^2C_{\Omega}t_n^\alpha}
\bra{\tau^{\min\{2-\alpha,\gamma\sigma\}}+t_n^{\alpha/2}h^2},
\end{align}
where $C_u$ may depend on $u$ and $T$, but is uniformly bounded with
respect to $\alpha$~and $\sigma$. It achieves an optimal time
accuracy of order $\mathcal{O}(\tau^{2-\alpha})$ if
$\gamma\geq\max\{1,(2-\alpha)/\sigma\}$.
\end{theorem}

\begin{remark}\label{remark: H2 analysis}
The special treatment of consistency error in time is motivated by the time-space error-splitting technique proposed originally in \cite{LiaoSunShi:2010,LiaoSunShi:2010B,LiaoSun:2010}
for obtaining the maximum norm error estimate via the discrete energy approach, see also \cite{LiaoYanZhang:2018}
for a recent application in the numerical analysis of a nonlinear subdiffusion problem.
To see it more clearly, we introduce $w=\brab{c(x)-\opL} u$ and reformulate the subdiffusion problem \eqref{eq: IBVP} into
\begin{equation*}
\begin{aligned}
w&=\fd{\alpha} u-f(x, t)&&\text{for $x\in\overline{\Omega}$ and $0<t\le T$,}\\
w&=\brab{c(x)-\opL} u&&\text{for $x\in\Omega$ and $0\le t\le T$.}
\end{aligned}
\end{equation*}
The fully discrete system follows as
\begin{equation*}
\begin{aligned}
    w_h^n&=(\dfd{\alpha}u_h)^{n}
        -f(x,t_{n})&
    &\text{for $x\in\overline{\Omega}_h$ and $1\le n\le N$,}\\
    w^n_h&=\brab{c(x)-\opL_h} u_h^{n}&
    &\text{for $x\in\Omega_h$ and $0\le n\le N$.}
\end{aligned}
\end{equation*}
Then the solution errors, $\tilde{w}^n_h:=w(x,t_n)-w^n_h$ and $\tilde{u}^n_h:=u(x,t_n)-u^n_h$
for~$x\in\overline{\Omega}_h$ satisfy
\begin{equation*}
\begin{aligned}
    \tilde{w}_h^n&=(\dfd{\alpha}\tilde{u}_h)^{n}
          -\overline{\Upsilon}_h^{n}[u]&
    &\text{for $x\in\overline{\Omega}_h$ and $1\le n\le N$,}\\
    \tilde{w}^n_h&=\brab{c(x)-\opL_h} \tilde{u}_h^{n}+R^n_s&
    &\text{for $x\in\Omega_h$ and $0\le n\le N$.}
\end{aligned}
\end{equation*}
We see that, the time and space truncation errors are redefined directly via this coupled error system. This is,
the time truncation error is defined uniformly over the closed space
domain and the spatial truncation error is defined uniformly over all time levels.

For the $H^1$-norm error estimate considered here, it needs only to redefine
the time consistency error as done in \eqref{eq: L1 time truncation}. It also motivates that
we can obtain an optimal $H^1$-norm error estimate via two stages: a time-discrete system
is considered in the first stage so that the time truncation error is defined uniformly
with respect to the spatial domain. As the spatial approximation of an elliptic problem,
the fully-discrete system can be treated traditionally in the second stage
and an optimal $H^1$-norm error estimate would be achieved because
it does not involve the time consistency error. We will illuminate the two-stage process
in the next section for a second-order scheme although it seems unusual in finite difference method.
\end{remark}

\section{Sharp $H^1$-norm error estimate for FracCN scheme}\label{sec:Alikhanov approach}
\setcounter{equation}{0}

To present an alternative approach
for a sharp $H^1$-norm error estimate, we recall the usual inner product
$\bra{v,w}=\int_{\Omega}v(x)w(x)\zd x$ with the associated $L^2(\Omega)$ norm $\mynorm{v}_{L^2(\Omega)}=\sqrt{\bra{v,v}}$.
For any functions $v$~and $w$ belonging to the space of grid functions that vanish
on the boundary~$\partial\Omega$,
define the $H^1$-seminorm $\abs{v}_{H^1(\Omega)}=\sqrt{\bra{\mu\partial_xv,\partial_xv}}\,.$
There exists a positive constant $C_{\Omega}$ is dependent
on the domain $\Omega$, the constants $\mu_0$ and $\mu_1$, such that
$\mynorm{v}_{L^2(\Omega)}\le C_{\Omega}\abs{v}_{H^1(\Omega)}$. Moreover, one has
\begin{align}\label{eq: continuous H1 norm}
\abs{v}_1\le C_{\Omega}\abs{v}_{H^1(\Omega)}
\end{align}
which can be checked by the Cauchy-Schwarz inequality with
$\partial_hv(x_{i-\frac12})=\frac1{h}\int_{x_{i-1}}^{x_i}v'(x)\zd x$.

\subsection{Nonuniform Alikhanov approximation}

Now we recall the nonuniform Alikhanov approximation,
see also \cite{LiaoMcLeanZhang:2018b}.
Let $\Pi_{1,k}v$ denote the linear interpolant of a function~$v$ with
respect to the nodes $t_{k-1}$~and $t_k$, and let $\Pi_{2,k}v$ denote the
quadratic interpolant with respect to $t_{k-1}$, $t_k$~and $t_{k+1}$.
It is easy to find that
\[
\bra{\Pi_{1,k}v}'(t)=\frac{\diff v^k}{\tau_k}
\;\;\text{and}\;\;
\bra{\Pi_{2,k}v}'(t)=\frac{\diff v^k}{\tau_k}
    +\frac{2(t-t_{k-1/2})}{\tau_k(\tau_k+\tau_{k+1})}
       \bra{\rho_k\diff v^{k+1}-\diff v^k}.
\]
The nonuniform Alikhanov formula to the Caputo derivative~$(\fd{\alpha}v)(t_{n-\theta})$ is defined by
\begin{align}\label{eq: Alikhanov}
(\dfd{\alpha}v)^{n-\theta}
    &\defeq\int_{t_{n-1}}^{t_{n-\theta}}\omega_{1-\alpha}(t_{n-\theta}-s)
        \bra{\Pi_{1,n}v}'(s)\zd{s}  +\sum_{k=1}^{n-1}\int_{t_{k-1}}^{t_k}   \omega_{1-\alpha}(t_{n-\theta}-s)\bra{\Pi_{2,k}v}'(s)\zd{s}
\nonumber\\
   &=a^{(n)}_0\diff v^n+\sum_{k=1}^{n-1}\braB{
    a^{(n)}_{n-k}\diff v^k+\rho_k b^{(n)}_{n-k}\diff v^{k+1}
        -b^{(n)}_{n-k}\diff v^k},
\end{align}
where the discrete coefficients
$a_{n-k}^{(n)}$ and $b_{n-k}^{(n)}$ are defined by
\begin{align}
&a^{(n)}_0\defeq\int_{t_{n-1}}^{t_{n-\theta}}
    \frac{\omega_{1-\alpha}(t_{n-\theta}-s)}{\tau_n}\zd{s}\;\;\mbox{and}\;\;
a^{(n)}_{n-k}\defeq\int_{t_{k-1}}^{t_k}\!\!\frac{\omega_{1-\alpha}(t_{n-\theta}-s)}{\tau_n}\zd{s}, \;\; 1\le k\leq n-1;\label{eq: an}\\
&b^{(n)}_{n-k}\defeq\frac{2}{\tau_k(\tau_k+\tau_{k+1})}\int_{t_{k-1}}^{t_k}
    (s-t_{k-\frac12})\omega_{1-\alpha}(t_{n-\theta}-s)\zd{s}, \quad 1\le k\leq n-1. \label{eq: bn}
\end{align}

Notice that while $\alpha\to1$, we have $\omega_{2-\alpha}(t)\to1$
and $\omega_{1-\alpha}(t)\to0$, uniformly for~$t$ in any compact
subinterval of the open half-line~$(0,\infty)$.  Thus, $a^{(n)}_0=\omega_{2-\alpha}((1-\theta)\tau_n)/\tau_n\to1/\tau_n$
whereas $a^{(n)}_{n-k}\to0$~and $b^{(n)}_{n-k}\to0$ for $1\le k\le n-1$.
It follows that $(\dfd{\alpha}v)^{n-\theta}\to\diff{v^n}/\tau_k$
and $\theta=\alpha/{2}\to1/2$
so the time-stepping scheme~\eqref{eq: discrete IBVP} with $\nu=\theta$ tends to the
classical second-order Crank--Nicolson method for a (classical) linear
reaction-diffusion equation. This is why we also call \eqref{eq: discrete IBVP} for the case $\nu=\theta$ as
a fractional Crank--Nicolson method.

Rearranging the terms in~\eqref{eq: Alikhanov}, we obtain
the compact form \eqref{eq: discrete Caputo} with $\nu=\theta$,
where the discrete convolution kernel $A_{n-k}^{(n,\theta)}$ is defined as follows:
$A_0^{(1,\theta)}\defeq a_0^{(1)}$ if $n=1$ and, for $n\geq2$,
\begin{equation}\label{eq: Alikhanov weights}
A^{(n,\theta)}_{n-k}\defeq\begin{cases}
    a^{(n)}_0+\rho_{n-1}b^{(n)}_1,
    &\text{for $k=n$,}\\
    a^{(n)}_{n-k}+\rho_{k-1}b^{(n)}_{n-k+1}-b^{(n)}_{n-k},
    &\text{for $2\le k\le n-1$,}\\
    a^{(n)}_{n-1}-b^{(n)}_{n-1},
    &\text{for $k=1$.}
\end{cases}
\end{equation}
Some useful properties of $A_{n-k}^{(n,\theta)}$ have been established recently by assuming that
\begin{itemize}
\item[\Ass{3r}.] The parameter $\theta=\alpha/2$, and the maximum time-step
ratio $\rho=7/4$.
\end{itemize}
\begin{theorem}\label{thm: FracCN kernels} \cite[Theorem 2.2]{LiaoMcLeanZhang:2018b}
If \Ass{3r} holds, then the discrete kernels in~\eqref{eq: Alikhanov weights} fulfills
\begin{itemize}
\item[(I)] The discrete kernels $A^{(n,\theta)}_{n-k}$ are positive and monotone,
\[A^{(n,\theta)}_{n-k-1}-A^{(n,\theta)}_{n-k}\ge
(1+\rho_k)b^{(n)}_{n-k}-\frac1{5\tau_k}\int_{t_{k-1}}^{t_k}\bra{t_k-s}\omega_{-\alpha}(t_{n-\theta}-s)\zd{s}>0
\quad\text{for $1\le k\le n-1;$}\]
\item[(II)] And, $ A^{(n,\theta)}_{0}-A^{(n,\theta)}_{1}>\theta\brab{2A^{(n,\theta)}_{0}-A^{(n,\theta)}_{1}}$ for~$n\ge2$;
\item[(III)] The discrete kernels $A^{(n,\theta)}_{n-k}$ are bounded, $A_{n-k}^{(n,\theta)}<A_{0}^{(n,\theta)}\leq
\frac{24}{11\tau_n}\int_{t_{n-1}}^{t_{n}}\omega_{1-\alpha}(t_n-s)\zd
s$ and
 $$A_{n-k}^{(n,\theta)}\geq\frac4{11\tau_k}\int_{t_{k-1}}^{t_{k}}\omega_{1-\alpha}(t_n-s)\zd s\quad\text{for $1\le k\le n.$}$$
\end{itemize}
\end{theorem}

The first two parts (I)-(II) ensures that \Ass{1} is valid, and the
last part (III) implies that \Ass{2} holds with
$\pi_A=\frac{11}{4}$. Hence we can use the complementary discrete
convolution kernel~$P^{(n,\theta)}_{n-k}$, see \eqref{eq: P
A}-\eqref{eq: P bound}, in this section. They allow us to apply
Lemma \ref{lem: general Dv v} and Theorem \ref{thm: H1 gronwall} and
establish the stability of the FracCN scheme~\eqref{eq: discrete
IBVP}. Actually, Theorems \ref{thm: FracCN kernels} and \ref{thm: stability}
imply the $H^1$-norm stability of the FracCN scheme for
the linear problem \eqref{eq: IBVP}.
\begin{corollary}\label{corol: FracCN stability}
If the local mesh restriction \Ass{3r} holds, then the FracCN method
\eqref{eq: discrete IBVP} with $\nu=\theta$ is unconditionally
stable in the discrete $H^1$-norm.
\end{corollary}

To derive a sharp $H^1$-norm error estimate, we need the following
two Lemmas.
\begin{lemma}\label{lem: Upsilon global singular}
Let $v\in\mathrm{C}_{\sigma}^{3}((0,T])$ with $\sigma\in(0,1)\cup(1,2)$.
If the mesh condition  \Ass{3r} holds,
then the local consistency error
$\Upsilon^{n-\theta}[v]$ of the nonuniform Alikhanov formula
$(\dfd{\alpha}v)^{n-\theta}$ in \eqref{eq: Alikhanov} with the discrete convolution kernels \eqref{eq: Alikhanov weights}
satisfies
\begin{align*}
\absb{\Upsilon^{n-\theta}[v]}\leq A_{0}^{(n,\theta)}\Gloc^n+
\sum_{k=1}^{n-1}\brab{A_{n-k-1}^{(n,\theta)}-A_{n-k}^{(n,\theta)}}\Ghis^k\quad \text{for $1\leq n\leq N$}
\end{align*}
where
\begin{align*}
\Gloc^k&\defeq \frac32\int_{t_{k-1}}^{t_{k-1/2}}(s-t_{k-1})^2|v'''(s)|\zd{s}
    +\frac{3\tau_k}{2}\int_{t_{k-1/2}}^{t_k}(t_k-s)|v'''(s)|\zd{s},\\
\Ghis^k&\defeq\frac52\int_{t_{k-1}}^{t_k}(s-t_{k-1})^2|v'''(s)|\zd{s}
    +\frac52\int_{t_k}^{t_{k+1}}(t_{k+1}-s)^2|v'''(s)|\zd{s}.
\end{align*}
Thus the global consistency error can be bounded by
\begin{align*}
\sum_{j=1}^nP_{n-j}^{(n,\theta)}\absb{\Upsilon^{j-\theta}[v]}
\le C_v\braB{\tau_1^{\sigma}/\sigma+t_1^{\sigma-3}\tau_2^3
+\frac1{1-\alpha}\max_{2\leq k\leq n}t_{k}^{\alpha}t_{k-1}^{\sigma-3}\tau_k^{3}/\tau_{k-1}^{\alpha}}.
\end{align*}
\end{lemma}
\begin{proof}
See Theorem 3.4 and Lemma 3.6 in \cite{LiaoMcLeanZhang:2018b}.
\end{proof}

Next Lemma suggests that the time weighted operator will not lead to
any loss of the temporal accuracy in the $H^1$-norm error analysis,
although the solution is non-smooth near $t=0$.
\begin{lemma}\label{lem:time-Weighted}
Let $v\in\mathrm{C}_{\sigma}^{2}((0,T])$ with $\sigma\in(0,1)\cup(1,2)$.
The truncation error of $v^{n-\theta}$ satisfies
\[
\absb{v(t_{1-\theta})-v^{1-\theta}}\leq C_v\,\tau_1^\sigma/\sigma
\quad\text{and}\quad
\absb{v(t_{j-\theta})-v^{j-\theta}}\leq C_vt_{j-1}^{\sigma-2}\tau_j^{2}
\quad\text{for $2\leq j\leq N$}
\]
such that
\begin{align*}
\max_{1\leq j\leq n}\big\{t_j^{\alpha/2}\absb{v(t_{j-\theta})-v^{j-\theta}}\big\}
\leq&\,C_v\tau_1^{\sigma+\alpha/2}/\sigma
+C_v\max_{2\leq k\leq n}t_k^{\alpha/2}t_{k-1}^{\sigma-2}\tau_k^{2}\quad \text{for $1\leq n\leq N.$}
\end{align*}
\end{lemma}
\begin{proof}The Taylor expansion with the integral remainder gives, see also \cite[Lemma~2.5]{LiaoZhaoTeng:2016},
\begin{align*}
v^{j-\theta}-v(t_{j-\theta})=\theta\int_{t_{j-1}}^{t_{j-\theta}}(s-t_{j-1})v''(s)\zd{s}
+(1-\theta)\int^{t_{j}}_{t_{j-\theta}}(t_{j}-s)v''(s)\zd{s}\,,\quad 1\leq j\leq N.
\end{align*}
The claimed results then follow immediately.
\end{proof}

\subsection{Two-stage convergence analysis}

Now we describe an alternative two-stage process for obtaining a sharp $H^1$-norm error estimate for
the second-order FracCN method \eqref{eq: discrete IBVP} with $\nu=\theta$ by assuming
that the subdiffusion problem \eqref{eq: IBVP} has a unique solution
$u\in \mathrm{C}\brab{[0,T];H^4(\Omega)}
\cap\mathrm{C}_{\sigma}^{2}\brab{(0,T];H^2(\Omega)}
\cap\mathrm{C}_{\sigma}^{3}\brab{(0,T];H^1(\Omega)}$.

\paragraph{Temporal error analysis via a time-discrete system}
We apply the nonuniform Alikhanov formula $(\dfd{\alpha}v)^{n-\theta}$ with the discrete
convolution kernels \eqref{eq: Alikhanov weights} to approximate the problem \eqref{eq: IBVP},
\begin{equation}\label{eq: time-discrete IBVP}
\begin{aligned}
(\dfd{\alpha}u)^{n-\nu}+\opL u^{n-\nu}
    &=c(x) u^{n-\nu}+f(x,t_{n-\nu})&
    &\text{for $x\in\Omega$ and $1\le n\le N$,}\\
    u^n&=u_b(x,t_n)&
    &\text{for $x\in\partial\Omega$ and $1\le n\le N$,}\\
    u^0&=u_0(x)&&\text{for $x\in\Omega$.}
\end{aligned}
\end{equation}
Then the solution error, $e^n=u(x,t_n)-u^n$ for $x\in \Omega$, satisfies
the zero-valued initial-boundary conditions and the governing equation
\begin{equation}\label{eq: FracCN time error IBVP}
(\dfd{\alpha}e)^{n-\theta}+\opL e^{n-\theta}
    =c(x)e^{n-\theta}+\Upsilon^{n-\theta}[u]+R_w^{n-\theta}\quad
    \text{for $x\in\Omega$ and $1\le n\le N$,}
\end{equation}
where $\Upsilon^{n-\theta}[u]$ is defined by \eqref{eq: Truncation discrete Caputo}
and $R^{n-\nu}_w\defeq\brab{c(x)-\opL}\kbrab{u^{n-\nu}-u(t_{n-\nu})}$ for $x\in \Omega$.

By taking the (continuous) inner product of the error equation
in~\eqref{eq: FracCN time error IBVP} with~$\opL e^{n-\theta}$, one applies the first Green formula to find
\begin{align*}
\brab{\bra{\dfd{\alpha}e}^{n-\theta},\opL e^{n-\theta}}
=&\,-\mynormb{\opL e^{n-\theta}}_{L^2(\Omega)}^2
+\brab{ce^{n-\theta},\opL e^{n-\theta}}
+\brab{R_w^{n-\theta},\opL e^{n-\theta}}\\
&\,+\brab{\mu\partial_x\Upsilon^{n-\theta}[u],\partial_xe^{n-\theta}}\\
\leq&\,\frac{1}2\kappa^2\mynormb{e^{n-\theta}}_{L^2(\Omega)}^2+\frac{1}2\mynormb{R_w^{n-\theta}}_{L^2(\Omega)}^2
+\absb{e^{n-\theta}}_{H^1(\Omega)}\absb{\Upsilon^{n-\theta}[u]}_{H^1(\Omega)}\\
\leq&\,
\frac{1}2\kappa^2C_{\Omega}\absb{e^{n-\theta}}_{H^1(\Omega)}^2+
\absb{e^{n-\theta}}_{H^1(\Omega)}\absb{\Upsilon^{n-\theta}[u]}_{H^1(\Omega)}+\frac12\mynormb{R_w^{n-\theta}}_{L^2(\Omega)}^2,
\end{align*}
where the Cauchy-Schwarz inequality and the embedding inequality have been used.
We apply Lemma \ref{lem: general Dv v} together with $\nu=\theta$ and Theorem \ref{thm: FracCN kernels} to obtain
\begin{align*}
\sum_{k=1}^nA^{(n,\theta)}_{n-k}\diff\brab{\absb{e^{k}}_{H^1(\Omega)}^2}
    \le&\,\kappa^2C_{\Omega}\absb{e^{n-\theta}}_{H^1(\Omega)}^2+
2\absb{e^{n-\theta}}_{H^1(\Omega)}\absb{\Upsilon^{n-\theta}[u]}_{H^1(\Omega)}+\mynormb{R_w^{n-\theta}}_{L^2(\Omega)}^2,
\end{align*}
which has the form of \eqref{eq: first Gronwall} with $\lambda_0:=\kappa^2C_{\Omega}$, $\lambda_l:=0$ $(l\ge1)$,
$$v^k:=\absb{e^{k}}_{H^1(\Omega)}, \quad\xi^n:=2\absb{\Upsilon^{n-\theta}[u]}_{H^1(\Omega)}\quad \text{and} \quad \eta^n:=\mynormb{R_w^{n-\theta}}_{L^2(\Omega)}.$$
Therefore, applying Theorem \ref{thm: improved gronwall} with $\pi_A=11/4$, we get
\begin{align*}
\absb{e^n}_{H^1(\Omega)}\le&\,4E_\alpha\brab{10\kappa^2C_{\Omega}t_n^\alpha}
    \max_{1\le k\le n}\sum_{j=1}^k P^{(k,\theta)}_{k-j}\absb{\Upsilon^{j-\theta}[u]}_{H^1(\Omega)}\\
  &\,+4\sqrt{\Gamma(1-\alpha)}E_\alpha\brab{10\kappa^2C_{\Omega}t_n^\alpha}
  \max_{1\le k\le n}t_k^{\alpha/2}\mynormb{R_w^{k-\theta}}_{L^2(\Omega)}\nonumber
\end{align*}
if the local assumption \Ass{3r} holds with the maximum time-step
size
$\taumax\le1/\sqrt[\alpha]{6\Gamma(2-\alpha)\kappa^2C_{\Omega}}\,$.
Then, applying Lemma \ref{lem: Upsilon global singular} (with
$v=\partial_xu$) and Lemma \ref{lem:time-Weighted}, one obtains
\begin{align}\label{eq: FracCN H1 time error}
\absb{u(t_n)-u^n}_{H^1(\Omega)}&\leq \frac{C_u}{\sigma(1-\alpha)}\braB{\tau_1^{\sigma}
+\max_{2\leq k\leq n}t_{k}^{\alpha}t_{k-1}^{\sigma-3}\tau_k^{3}/\tau_{k-1}^{\alpha}
    +\max_{2\leq k\leq n}t_k^{\alpha/2}t_{k-1}^{\sigma-2}\tau_k^{2}}.
\end{align}

\paragraph{Spatial error analysis via the fully-discrete system} Now return to
the fully-discrete system \eqref{eq: discrete IBVP} with $\nu=\theta$, which can be viewed as
the spatial approximation of time-discrete system  \eqref{eq: time-discrete IBVP}.
Under our priori assumptions to the problem \eqref{eq: IBVP},
this system has a unique solution $u^n\in H^4(\Omega)$ for $1\leq n\leq N$.
Thus the solution error, $z^n_h:=u^n-u_h^n$ for~$x\in\overline{\Omega}_h$,
satisfies the zero-valued initial-boundary conditions, and the governing equation
\begin{equation}\label{eq: FracCN space error IBVP}
(\dfd{\alpha}z_h)^{n-\theta}+\opL_hz_h^{n-\theta}
    =c(x)z_h^{n-\theta}+R_s^{n-\theta}\quad
    \text{for $x\in\Omega_h$ and $1\le n\le N$,}
\end{equation}
where $R_s^{n-\theta}$ is defined by \eqref{eq: space truncation}.
We will proceed to apply the standard $H^1$-norm analysis, as done in the subsection 2.2.
By taking the inner product of~\eqref{eq: FracCN space error IBVP} with~$\bra{\dfd{\alpha}z_h}^{n-\theta}$, one has
\begin{align*}
\myinnerb{\opL_hz_h^{n-\theta},\bra{\dfd{\alpha}z_h}^{n-\theta}}
&=-\mynormb{\bra{\dfd{\alpha}z_h}^{n-\theta}}^2+
\myinnerb{cz_h^{n-\theta},\bra{\dfd{\alpha}z_h}^{n-\theta}}+\myinnerb{R_s^{n-\theta},\bra{\dfd{\alpha}z_h}^{n-\theta}}\\
&\leq\frac12\kappa^2\mynormb{z_h^{n-\theta}}^2+\frac12\mynormb{R_s^{n-\theta}}^2\,.
\end{align*}
Therefore, applying Lemma \ref{lem: general Dv v} and the embedding inequality, one gets
\begin{align*}
\sum_{k=1}^nA^{(n,\theta)}_{n-k}\diff\brab{\absb{z_h^k}_1^2}
    \le\kappa^2C_{\Omega}\braB{(1-\theta)\absb{z_h^n}_1+\theta\absb{z_h^{n-1}}_1}^2
    +\mynormb{R_s^{n-\theta}}^2\quad\text{for $1\le n\le N$},
\end{align*}
which has the form of \eqref{eq: first Gronwall} with
$\lambda_0:=\kappa^2C_{\Omega}$, $\lambda_l:=0$ $(l\ge1)$,
$v^k:=\absb{z_h^{k}}_{1}$, $\xi^n:=0$ and
$\eta^n:=\mynormb{R_s^{n-\theta}}.$
Then the fractional Gr\"{o}nwall inequality in Theorem \ref{thm:
improved gronwall} (taking $\rho=7/4$ and $\pi_A=11/4$) and the
error estimate $\mynormb{R_s^{k-\theta}}\le C_uh^2$  yield
\begin{align}\label{eq: FracCN H1 space error}
\absb{u^n-u_h^n}_1\le&\,4E_\alpha\brab{10\kappa^2C_{\Omega}t_n^\alpha}
    \sqrt{\Gamma(1-\alpha)}\max_{1\le k\le n}\big\{t_k^{\alpha/2}\mynormb{R_s^{k-\theta}}\big\}\\
    \le&\, C_uE_\alpha\brab{10\kappa^2C_{\Omega}t_n^\alpha}\sqrt{\Gamma(1-\alpha)}t_n^{\alpha/2}h^2,\nonumber
\end{align}
if the assumption \Ass{3r} holds with the maximum time-step size
$\taumax\le1/\sqrt[\alpha]{6\Gamma(2-\alpha)\kappa^2C_{\Omega}}\,$.

We are in the position to complete the error estimate.
Combining \eqref{eq: FracCN H1 time error} with \eqref{eq: FracCN H1 space error},
one can apply the triangle inequality and the relationship \eqref{eq: continuous H1 norm} to find
\begin{align}\label{eq: FracCN H1 error}
\absb{u(t_n)-u_h^n}_1\le&\,\absb{u(t_n)-u^n}_1+\absb{u^n-u_h^n}_1
\leq C_{\Omega}\absb{u(t_n)-u^n}_{H^1(\Omega)}+\absb{u^n-u_h^n}_1\\
\le&\,\frac{C_u}{\sigma(1-\alpha)}\braB{\tau_1^{\sigma}
+\max_{2\leq k\leq n}t_{k}^{\alpha}t_{k-1}^{\sigma-3}\tau_k^{3}/\tau_{k-1}^{\alpha}
    +\max_{2\leq k\leq n}t_k^{\alpha/2}t_{k-1}^{\sigma-2}\tau_k^{2}+t_n^{\alpha/2}h^2},\nonumber
\end{align}
where $C_u$ may depend on $u$ and $T$, but is uniformly bounded with respect to $\alpha$~and $\sigma$.
If the mesh assumption \Mss{-conv} holds, then
$\tau_1\le C_{\gamma}\tau^{\gamma}$ and
\begin{align}\label{eq: D tau order}
t_{k}^{\alpha}t_{k-1}^{\sigma-3}\tau_k^{3}/\tau_{k-1}^{\alpha}
&=t_{k}^{2\alpha}t_{k-1}^{\sigma-3-\alpha}\tau_k^{3-\alpha}(\tau_k/t_k)^{\alpha}(\tau_{k-1}/t_{k-1})^{-\alpha}
\leq C_{\gamma}t_{k}^{2\alpha}t_{k-1}^{\sigma-3-\alpha}\tau_k^{3-\alpha}\\
&\leq C_{\gamma}t_{k}^{\sigma-3+\alpha}\tau_k^{3-\alpha}\leq C_\gamma t_k^{\sigma-3+\alpha}\tau_k^{3-\alpha-\beta}
    \bigl(\tau\min\{1,t_k^{1-1/\gamma}\}\bigr)^\beta\nonumber\\
&\leq C_\gamma t_k^{\sigma-\beta/\gamma}(\tau_k/t_k)^{3-\alpha-\beta}\tau^\beta
\leq C_\gamma t_k^{\max\{0,\sigma-(3-\alpha)/\gamma\}}\tau^{\beta},
    \quad 2\leq k\leq n; \nonumber
\end{align}
where~$\beta:=\min\{2,\gamma\sigma\}$. In addition,
\begin{align}\label{eq: offset order}
t_k^{\alpha/2}t_{k-1}^{\sigma-2}\tau_k^2
    &\leq C_\gamma t_k^{\sigma-2+\alpha/2}\tau_k^{2-\beta}
        \bigl(\tau\min\{1,t_k^{1-1/\gamma}\}\bigr)^\beta\\
&\le C_\gamma t_k^{\sigma+\alpha/2-\beta/\gamma}(\tau_k/t_k)^{2-\beta}\tau^\beta
\le C_\gamma t_k^{\alpha/2+\max\{0,\sigma-2/\gamma\}}\tau^\beta,\quad 2\leq k\le n.\nonumber
\end{align}
So the following result is achieved by inserting \eqref{eq: D tau order} and \eqref{eq: offset order} into \eqref{eq: FracCN H1 error}.

\begin{theorem}\label{thm: FracCN convergence}
Suppose that the initial-boundary value problem \eqref{eq: IBVP} of
the subdiffusion equation has a solution $u\in
\mathrm{C}\brab{[0,T];H^4(\Omega)}
\cap\mathrm{C}_{\sigma}^{2}\brab{(0,T];H^2(\Omega)}
\cap\mathrm{C}_{\sigma}^{3}\brab{(0,T];H^1(\Omega)}$, and consider
the fractional Crank-Nicoslon method \eqref{eq: discrete IBVP} using
the Alikhanov formula $(\dfd{\alpha}v)^{n-\theta}$ with the discrete
convolution kernels \eqref{eq: Alikhanov weights}. If the local mesh
condition \Ass{3r} holds with the maximum time-step size
$\taumax\le1/\sqrt[\alpha]{6\Gamma(2-\alpha)\kappa^2C_{\Omega}}$,
then the discrete solution $u_h^n$ is convergent in the
discrete $H^1$-norm,
\[
\absb{u(t_n)-u_h^n}_1
    \le \frac{C_u}{\sigma(1-\alpha)}\braB{\tau_1^{\sigma}
+\max_{2\leq k\leq n}t_{k}^{\alpha}t_{k-1}^{\sigma-3}\tau_k^{3}/\tau_{k-1}^{\alpha}
    +\max_{2\leq k\leq n}t_k^{\alpha/2}t_{k-1}^{\sigma-2}\tau_k^{2}+t_n^{\alpha/2}h^2}.
\]
In particular, if the mesh assumption \Mss{-conv} holds, then
\[
\absb{u(t_n)-u_h^n}_1\le
\frac{C_u}{\sigma(1-\alpha)}\brab{\tau^{\min\{\gamma\sigma,2\}}+h^2}
    \quad\text{for $1\le n\le N$,}
\]
where $C_u$ may depend on $u$ and $T$, but is uniformly bounded with respect to $\alpha$~and $\sigma$.
\end{theorem}

\begin{remark}\label{remark:Alikhanov-H1consistence}
As noted early in \cite{LiaoMcLeanZhang:2018b}, by an argument
similar to that in \eqref{eq: D tau order}, it is not difficult to
show that
$t_{k}^{\alpha}t_{k-1}^{\sigma-3}\tau_k^{3}/\tau_{k-1}^{\alpha} \le
C_{\gamma}t_k^{\sigma-(3-\alpha)/\gamma}\tau^{3-\alpha}$, which
means that the Alikhanov formula $(\dfd{\alpha}v)^{n-\theta}$
approximates $(\fd{\alpha}u)(t_{n-\theta})$ to
order~$\mathcal{O}(\tau^{3-\alpha})$ if
$\gamma\ge(3-\alpha)/\sigma$.  However, the term~\eqref{eq: offset
order} arising from~$R_w^{n-\theta}$ would still limit the
convergence rate for the overall scheme to order~$O(\tau^{2})$.

From the point of view of different spatial discretization methods,
the two-stage analysis would be more general that
the direct error splitting technique in subsection 3.2.
On the other hand, the traditional $H^1$-norm analysis in subsection 3.1 will
yield a suboptimal error estimate because the global consistency error
$\sum_{j=1}^k P^{(k,\theta)}_{k-j}\mynormb{\Upsilon_h^{j-\theta}[u]}^2$ also has a loss of time accuracy.
Actually, by using the discrete convolution bound of the local consistence error in Lemma \ref{lem: Upsilon global singular},
one can present an proof similar to that of Lemma \ref{lemma:L1formulaNonuniform-H1consistence}
and find the following estimate.
\begin{lemma}\label{lemma:Alikhanov-H1consistence}
 If $v\in \mathrm{C}_{\sigma}^{3}((0,T])$ for $\sigma\in(\frac{\alpha}2,1)\cup(1,2)$ and the maximum step ratio $\rho\leq1$, then
 the global consistency error of the nonuniform Alikhanov formula
$(\dfd{\alpha}v)^{n-\theta}$ in \eqref{eq: Alikhanov} with the discrete convolution kernels \eqref{eq: Alikhanov weights} satisifies
\begin{align*}
\sum_{j=1}^nP^{(n,\theta)}_{n-j}\absb{\Upsilon^{j-\theta}[v]}^2\leq
\frac{C_v}{\sigma^2}\tau_1^{2\sigma-\alpha}+t_1^{2\sigma-6-\alpha}\tau_2^6
+\frac{C_v}{1-\alpha}\max_{2\leq k\leq n}t_k^{\alpha}t_{k-1}^{2\sigma-6}\tau_k^{6}/\tau_{k-1}^{2\alpha}
\quad\text{for~$1\le n\le N$.}
\end{align*}
Moreover, if the time mesh satisfies \Mss{-conv}, then
\begin{align*}
\sum_{j=1}^nP^{(n,\theta)}_{n-j}\absb{\Upsilon^{j-\theta}[v]}^2\leq
\frac{C_v}{\sigma^2(1-\alpha)}\tau^{2\min\{2,\gamma(\sigma-\alpha/2)\}}\quad\text{for~$1\le n\le N$.}
\end{align*}
\end{lemma}
\end{remark}

\section{Numerical examples}\label{sec:numerical}
\setcounter{equation}{0}

We present some numerical results to verify our error estimates. Always,
consider the reaction-subdiffusion problem \eqref{eq: IBVP} in the
spatial domain $\Omega=(0,\pi)$ and the time interval $[0,T]$ with $T=1$.
In the computations, the domain $(0,\pi)$ is divided into $M$ equally
spaced subintervals with a mesh length $h=\pi/M$, and the time
interval $[0,1]$ is divided into $N$ parts by an initially graded grid \eqref{eq: graded mesh}
with $T_0=\min\{\gamma^{-1},2^{-\gamma}\}$. Throughout our tests, we measure the discrete $H^1$-seminorm solution error
$e(M,N)=\max_{1\leq n\leq N}|u(t_n)-u_h^n|_1$.
Since the convergence behavior of the spatial discretization is well
understood, we focus on the temporal convergence here by setting a sufficiently large $M$
such that the time error dominates the spatial error in each run and $e(M,N)\approx e(N)$.
The experimental rate (list as ``Order" in tables) in temporal direction is
estimated by using
$\hbox{Order}=\log_2\bra{e(N)/e(2N)}.$

{\bf Example 1.} {\bf Numerical results for the fully discrete L1
scheme.}  We set a diffusive coefficient
$\mu(x)=\exp(x)$,
a reaction coefficient $c(x)=2\sin(x)+1$, and a specific
source term $f(x,t)$
such that the exact solution $u(x,t)=\omega_{1+\sigma}(t)\sin(x)$.
It is seen that this solution
fulfills the assumption
$u\in\mathrm{C}_{\sigma}^{2}((0,T];H^1(\Omega))$ for the regularity
parameter $\sigma\in (0,1)\cup (1,2)$.


\begin{table}[ht!]
 \centering \caption{Numerical accuracy for Example 1 with $\sigma=2-\alpha$ and $\gamma=1$.}
\begin{tabular}[b]{ccccccccccccc}
\hline
{$N$}  & \multicolumn{2}{c}{$\alpha=0.1,\,\sigma=1.9$} &     & \multicolumn{2}{c}{$\alpha=0.5,\,\sigma=1.5$}  &  &\multicolumn{2}{c}{$\alpha=0.9,\,\sigma=1.1$}\\
\cline{2-3} \cline{5-6} \cline{8-9}
                    &$e(M,N)$              &   $\text{Order}$           &$\,$     &$e(M,N)$        &$\text{Order}$   &$\,$         &$e(M,N)$  &$\text{Order}$\\
\hline
      $ 100  $   &     3.84e-06  &  1.83       &$\,$     &  1.71e-04   &      1.38  &$\,$     &    1.03e-03 &        0.94   \\
      $ 200  $   &     1.08e-06  &  1.84       &$\,$     &  6.56e-05   &      1.40  &$\,$     &    5.36e-04 &        0.96   \\
      $ 400  $   &     3.02e-07  &  1.84       &$\,$     &  2.48e-05   &      1.42  &$\,$     &    2.75e-04 &        0.98   \\
      $ 800  $  &      8.46e-08  &  1.84       &$\,$     &  9.27e-06   &      1.43  &$\,$     &    1.40e-04 &        0.99   \\
      $ 1600  $  &     2.37e-08  &   *         &$\,$     &  3.43e-06   &  *         &$\,$     &    7.04e-05 &     *        \\
 \hline
 $\min\{\gamma\sigma,2-\alpha\}$    &   & 1.90     &          &              &         1.50    &        &               &    1.10  \\
  \hline
\end{tabular}\label{tablep:2}
\end{table}

 To test the sharpness of our error
estimate Theorem \ref{th:Convergence-nonuniformL1Scheme}, we
consider four different scenarios, respectively, in Tables
\ref{tablep:2}-\ref{tablep:5}. Setting the fixed and sufficiently
big $M=20000$, the sufficiently small value of $h$ can guarantee
that the dominated errors arise from the L1 approximation of Caputo
derivative. By taking $\sigma=2-\alpha$ and $\gamma=1$, the
computational results of the scheme for different
$\alpha=0.1,0.5,0.9$ are presented in Table \ref{tablep:2}.  It is
observed that the scheme has the temporal order
$\mathcal{O}(\tau^{2-\alpha})$, which is consistent with our theoretical
analysis.


\begin{table}[ht!]
 \centering \caption{Numerical accuracy for Example 1 with $\alpha=0.5$, $\sigma=0.5$ and $\gamma_{\text{opt}}=3$.}
\begin{tabular}[b]{ccccccccccccc}
\hline
{$N$}  & \multicolumn{2}{c}{$\gamma=1$} &     & \multicolumn{2}{c}{$\gamma=3$}&     & \multicolumn{2}{c}{$\gamma=3.75$}\\
\cline{2-3} \cline{5-6} \cline{8-9}
                    &$e(M,N)$              &   $\text{Order}$           &$\,$     &$e(M,N)$        &$\text{Order}$   &$\,$         &$e(M,N)$  &$\text{Order}$\\
\hline
      $ 100  $   &     2.57e-02  &     0.45   &$\,$     &    6.34e-04 &      1.43   &$\,$     &    4.66e-04 &        1.47   \\
      $ 200  $   &     1.88e-02  &     0.46   &$\,$     &    2.34e-04 &      1.45   &$\,$     &    1.68e-04 &        1.48   \\
      $ 400  $   &     1.37e-02  &     0.47   &$\,$     &    8.56e-05 &      1.47   &$\,$     &    6.01e-05 &        1.49   \\
      $ 800  $   &     9.88e-03  &     0.47   &$\,$     &    3.10e-05 &      1.48   &$\,$     &    2.14e-05 &        1.49  \\
      $ 1600  $  &     7.11e-03  &   *        &$\,$     &    1.11e-05 &  *          &$\,$     &    7.63e-06 &     *        \\
 \hline
 $\min\{\gamma\sigma,2-\alpha\}$    &   & 0.50     &          &              &         1.50    &        &               &    1.50  \\
  \hline
\end{tabular}\label{tablep:3}
\end{table}

\begin{table}[ht!]
 \centering \caption{Numerical accuracy for Example 1 with $\alpha=0.5$, $\sigma=0.75$ and $\gamma_{\text{opt}}=2$.}
\begin{tabular}[b]{ccccccccccccc}
\hline
{$N$}  & \multicolumn{2}{c}{$\gamma=1$} &     & \multicolumn{2}{c}{$\gamma=2$}&     & \multicolumn{2}{c}{$\gamma=2.5$}\\
\cline{2-3} \cline{5-6} \cline{8-9}
                    &$e(M,N)$              &   $\text{Order}$           &$\,$     &$e(M,N)$        &$\text{Order}$   &$\,$         &$e(M,N)$  &$\text{Order}$\\
\hline
      $ 100  $   &      3.70e-03  &      0.70    &$\,$     &     2.26e-04 &      1.41  &$\,$     &    1.47e-04 &        1.47  \\
      $ 200  $   &      2.28e-03  &      0.71    &$\,$     &     8.48e-05 &      1.43 &$\,$      &    5.30e-05 &        1.48   \\
      $ 400  $   &      1.39e-03  &      0.72    &$\,$     &     3.14e-05 &      1.45 &$\,$      &    1.90e-05 &        1.48  \\
      $ 800  $   &      8.46e-04  &      0.72    &$\,$     &     1.15e-05 &      1.46 &$\,$      &    6.79e-06 &        1.49 \\
      $ 1600  $  &      5.12e-04  &   *          &$\,$     &     4.18e-06 &  *        &$\,$      &    2.42e-06 &     *        \\
 \hline
 $\min\{\gamma\sigma,2-\alpha\}$    &   & 0.75     &          &              &         1.50   &        &               &    1.50  \\
  \hline
\end{tabular}\label{tablep:4}
\end{table}

\begin{table}[ht!]
 \centering \caption{Numerical accuracy for Example 1 with $\alpha=0.5$, $\sigma=1.25$ and $\gamma_{\text{opt}}=1.2$.}
\begin{tabular}[b]{ccccccccccccc}
\hline
{$N$}  & \multicolumn{2}{c}{$\gamma=1$} &     & \multicolumn{2}{c}{$\gamma=1.2$}&     & \multicolumn{2}{c}{$\gamma=1.8$}\\
\cline{2-3} \cline{5-6} \cline{8-9}
                    &$e(M,N)$              &   $\text{Order}$           &$\,$     &$e(M,N)$        &$\text{Order}$   &$\,$         &$e(M,N)$  &$\text{Order}$\\
\hline
      $ 100  $   &      2.75e-04  &      1.17   &$\,$     &    1.18e-04 &      1.39  &$\,$     &    7.76e-05 &        1.49    \\
      $ 200  $   &      1.22e-04  &      1.20   &$\,$     &    4.52e-05 &      1.41  &$\,$     &    2.75e-05 &        1.53    \\
      $ 400  $   &      5.33e-05  &      1.21   &$\,$     &    1.70e-05 &      1.43  &$\,$     &    9.55e-06 &        1.60    \\
      $ 800  $   &      2.31e-05  &      1.22   &$\,$     &    6.34e-06 &      1.44  &$\,$     &    3.14e-06 &        1.66   \\
      $ 1600  $  &      9.92e-06  &   *         &$\,$     &    2.34e-06 &  *        &$\,$      &    9.96e-07 &     *        \\
 \hline
 $\min\{\gamma\sigma,2-\alpha\}$    &   & 1.25     &          &              &         1.50    &        &               &    1.50  \\
  \hline
\end{tabular}\label{tablep:5}
\end{table}

Numerical results in Tables \ref{tablep:3}-\ref{tablep:5} (with
$\alpha=0.5$ and $\sigma<2-\alpha$) support the predicted time
accuracy in Theorem \ref{th:Convergence-nonuniformL1Scheme}. In the case of uniform
mesh $\gamma=1$, the solution is accurate of order
$\mathcal{O}(\tau^{\sigma})$, and nonuniform meshes improve the
numerical precision and convergence rate of solution. When
the grid parameter $\gamma\geq\gamma_{\text{opt}}$, the optimal time
accuracy $\mathcal{O}(\tau^{2-\alpha})$ is observed. Thus the $H^1$-norm error
estimate \eqref{L1convegenceerror} is sharp.

{\bf Example 2.} {\bf Numerical results for the fully discrete FracCN scheme.}
We choose  $\mu(x)=\cos(x)+2$,  $c(x)=2\sin(x)+1$, $u^0=\sin(x)$,
and a forcing source $f(x,t)$ such that the problem has a solution $ u(x,t)=(1+\omega_{1+\sigma}(t))\sin(x).$

\begin{table}[ht!]
 \centering \caption{Numerical accuracy for Example 2 with $\sigma=1+\alpha$ and $\gamma=1$.}
\begin{tabular}[b]{ccccccccccccc}
\hline
{$N$}  & \multicolumn{2}{c}{$\alpha=0.4,\,\sigma=1.4$} &     & \multicolumn{2}{c}{$\alpha=0.6,\,\sigma=1.6$}  &  &\multicolumn{2}{c}{$\alpha=0.8,\,\sigma=1.8$}\\
\cline{2-3} \cline{5-6} \cline{8-9}
                    &$e(M,N)$              &   $\text{Order}$           &$\,$     &$e(M,N)$        &$\text{Order}$   &$\,$         &$e(M,N)$  &$\text{Order}$\\
\hline
      $ 128  $   &   3.42e-05    &  1.63  &$\,$     &    3.43e-05 &  1.97  &$\,$     &    2.65e-05 &   1.97  \\
      $ 256  $   &   1.10e-05    &  1.57  &$\,$     &    8.73e-06 &  1.96  &$\,$     &    6.76e-06 &   1.96  \\
      $ 512  $   &   3.73e-06    &  1.54  &$\,$     &    2.23e-06 &  1.90  &$\,$     &    1.73e-06 &   1.92  \\
      $ 1024 $   &   1.28e-06    &  1.51  &$\,$     &    5.40e-07 &  1.71  &$\,$     &    4.61e-07 &   1.86  \\
      $ 2048 $   &   4.50e-07    &  1.49  &$\,$     &    1.33e-07 &  1.71  &$\,$     &    1.27e-07 &   1.83  \\
      $ 4096 $   &   1.60e-07    &   *    &$\,$     &    5.02e-08 &  *     &$\,$     &    3.57e-08 &     *        \\
 \hline
 $\min\{\gamma\sigma,2\}$     &   & 1.40     &          &              &         1.60    &        &               &    1.80  \\
  \hline
\end{tabular}\label{tablep:6}
\end{table}

\begin{table}[ht!]
 \centering \caption{Numerical accuracy for Example 2 with $\sigma=1.2$ and $\alpha=0.4$.}
\begin{tabular}[b]{ccccccccccccc}
\hline
{$N$}  & \multicolumn{2}{c}{$\gamma=1$} &     & \multicolumn{2}{c}{$\gamma=5/3=\gamma_{\text{opt}}$}  &  &\multicolumn{2}{c}{$\gamma=2$}\\
\cline{2-3} \cline{5-6} \cline{8-9}
                    &$e(M,N)$              &   $\text{Order}$           &$\,$     &$e(M,N)$        &$\text{Order}$   &$\,$         &$e(M,N)$  &$\text{Order}$\\
\hline
      $ 128  $   &     6.17e-05  &  1.36  &$\,$     &    1.32e-05 &   2.05  &$\,$     &   1.39e-05 &   2.04  \\
      $ 256  $   &     2.40e-05  &  1.34  &$\,$     &    3.19e-06 &   2.00  &$\,$     &   3.36e-06 &   2.04  \\
      $ 512  $   &     9.49e-06  &  1.31  &$\,$     &    7.98e-07 &   2.06  &$\,$     &   8.16e-07 &   2.06  \\
      $ 1024 $   &     3.83e-06  &  1.29  &$\,$     &    1.91e-07 &   2.05  &$\,$     &   1.96e-07 &   2.06  \\
      $ 2048 $   &     1.57e-06  &   *    &$\,$     &    4.61e-08 &  *      &$\,$     &   4.70e-08 &     *        \\
 \hline
 $\min\{\gamma\sigma,2\}$     &   & 1.20     &          &              &         2.00    &        &               &    2.00  \\
  \hline
\end{tabular}\label{tablep:7}
\end{table}

\begin{table}[ht!]
 \centering \caption{Numerical accuracy for Example 2 with $\sigma=0.8$ and $\alpha=0.4$.}
\begin{tabular}[b]{ccccccccccccc}
\hline
{$N$}  & \multicolumn{2}{c}{$\gamma=2$} &     & \multicolumn{2}{c}{$\gamma=5/2=\gamma_{\text{opt}}$}  &  &\multicolumn{2}{c}{$\gamma=3$}\\
\cline{2-3} \cline{5-6} \cline{8-9}
                    &$e(M,N)$              &   $\text{Order}$           &$\,$     &$e(M,N)$        &$\text{Order}$   &$\,$         &$e(M,N)$  &$\text{Order}$\\
\hline
      $ 128  $   &     2.43e-05  &  2.12  &$\,$     &    2.49e-05 &  2.12  &$\,$     &    2.75e-05  &    2.12 \\
      $ 256  $   &     5.59e-06  &  1.69  &$\,$     &    5.72e-06 &  2.15  &$\,$     &    6.35e-06  &    2.15 \\
      $ 512  $   &     1.74e-06  &  1.61  &$\,$     &    1.29e-06 &  2.35  &$\,$     &    1.43e-06  &    2.33 \\
      $ 1024 $   &     5.69e-07  &  1.61  &$\,$     &    2.53e-07 &  2.43  &$\,$     &    2.84e-07  &    2.33 \\
      $ 2048 $   &     1.87e-07  &   *    &$\,$     &    4.67e-08 &  *     &$\,$     &    5.66e-08  &     *        \\
 \hline
 $\min\{\gamma\sigma,2\}$     &   & 1.60     &          &              &         2.00    &        &               &    2.00  \\
  \hline
\end{tabular}\label{tablep:8}
\end{table}

\begin{table}[ht!]
 \centering \caption{Numerical accuracy for Example 2 with $\sigma=0.4$ and $\alpha=0.4$.}
\begin{tabular}[b]{ccccccccccccc}
\hline
{$N$}  & \multicolumn{2}{c}{$\gamma=2$} &     & \multicolumn{2}{c}{$\gamma=5/2$}  &  &\multicolumn{2}{c}{$\gamma=5=\gamma_{\text{opt}}$}\\
\cline{2-3} \cline{5-6} \cline{8-9}
                    &$e(M,N)$              &   $\text{Order}$           &$\,$     &$e(M,N)$        &$\text{Order}$   &$\,$         &$e(M,N)$  &$\text{Order}$\\
\hline
      $ 128  $   &      3.35e-03  &  0.81  &$\,$     &    1.50e-03&  1.01  &$\,$     &   4.56e-04 &   2.17  \\
      $ 256  $   &      1.91e-03  &  0.81  &$\,$     &    7.41e-04&  1.00  &$\,$     &   1.01e-04 &   2.20  \\
      $ 512  $   &      1.09e-03  &  0.81  &$\,$     &    3.71e-04&  1.00  &$\,$     &   2.20e-05 &   2.17  \\
      $ 1024 $   &      6.21e-04  &  0.80  &$\,$     &    1.85e-04&  1.00  &$\,$     &   4.90e-06 &   2.14  \\
      $ 2048 $   &      3.56e-04  &   *    &$\,$     &    9.25e-05&  *     &$\,$     &   1.11e-06 &     *        \\
 \hline
 $\min\{\gamma\sigma,2\}$     &   & 0.80     &          &              &        1.00    &        &               &    2.00  \\
  \hline
\end{tabular}\label{tablep:9}
\end{table}

The solution is approximated by the FracCN scheme \eqref{eq: discrete IBVP} with $\nu=\alpha/2$.
For different fractional order $\alpha$,
the numerical results are computed with varying temporal stepsizes
and fixed sufficiently large spatial points $M=20000$. Like before,
for fixed $M$, the computational errors and numerical convergence
orders in the $H^1$-norm are given in Tables
\ref{tablep:6}-\ref{tablep:9} with different temporal step sizes,
from which, the $\mathcal{O}\brab{\tau^{\min\{\gamma\sigma,2\}}}$
convergence of the difference scheme \eqref{eq: discrete IBVP} is
apparent, indicating the sharpness of our estimate in Theorem
\ref{thm: FracCN convergence}.

\end{document}